\RequirePackage{fix-cm}
\documentclass[smallextended]{svjour3}       
\smartqed  
\usepackage[dvipdfmx]{graphicx}

\usepackage{mathtools}
\usepackage[dvipdfmx]{color}
\usepackage{amsmath,amssymb,mathrsfs}
\usepackage{float}
\usepackage{bm}
\usepackage{subfig}
\mathtoolsset{showonlyrefs=true}

\newcommand{\p}{\partial}
\newcommand{\bs}[1]{\boldsymbol{#1}}
\newcommand{\VN}{\bs{V}_{\mathrm{Newton}}}
\newcommand{\hatVN}{\hat{\bs{V}}_{\mathrm{Newton}}}
\begin{document}

\title{Convergence of a finite difference scheme for the Kuramoto--Sivashinsky equation defined on an expanding circle
  \thanks{The authors are grateful to Professors Hiroshi Kokubu (Kyoto University) and Tomoyuki Miyaji (Kyoto University) for continuous discussions and their valuable advice in the development of this paper.
    This study is partially supported by JSPS KAKENHI,
Grant Numbers 20K22307 (S. Kobayashi) and 19H01807 (S. Yazaki).}
}

\titlerunning{Finite difference scheme for the KS equation defined on an expanding circle}

\author{Shunsuke Kobayashi \and Shigetoshi Yazaki.
}


\institute{S. Kobayashi \at
Center for Science Adventure and Collaborative Research Advancement, Graduate School of Science, Kyoto University, Kitashirakawa Oiwake-cho, Sakyo-ku, Kyoto 606-8502, Japan;
RIKEN iTHEMS, 2-1 Hirosawa, Wako-shi, Saitama 351-0198, Japan\\
  \email{kobayashi.shunsuke.6e@kyoto-u.ac.jp}           
  \and
  S. Yazaki \at
  Department of Mathematics, School of Science and Technology, Meiji University, 1-1-1 Higashi-Mita, Tama-ku, Kawasaki-shi, Kanagawa 214-8571, Japan
}

\date{Received: date / Accepted: date}

\maketitle

\begin{abstract}
  This paper presents a finite difference method combined with the Crank--Nicolson scheme of the Kuramoto--Sivashinsky equation defined on an expanding circle (\cite{KUY}), and the existence, uniqueness, and second-order error estimate of the scheme. The equation can be obtained as a perturbation equation from the circle solution to an interfacial equation and can provide guidelines for understanding the wavenumber selection of solutions to the interfacial equation. Our proposed numerical scheme can help with such a mathematical analysis.
  
\end{abstract}


\section{Introduction}\label{intro}
A mathematical analysis of the behavior of a gaseous combustion flame front has long been under study, starting with the pioneering research by Sivashinsky \cite{S}.
In recent years, the behavior of the flame front during the smoldering combustion on a thin solid (e.g. a sheet of paper) has been vigorously studied both experimentally (\cite{GKKY,KSTK,OBK,ZOM,ZM1,ZM2,ZRRG}) and mathematically (\cite{FMP,GKUY,GKY,IIMO1,IIMO2,IM1,IM2,KKUYB,KS}) aspects.
The Kuramoto--Sivashinsky (KS) equation (\cite{KT,S}), which is a well-known mathematical model of gaseous combustion, has been applied to research on smoldering combustion (\cite{GKKY,GKUY,GKY,KKUYB}).

The purpose of the present paper is to numerically solve the following time evolution equation (\cite{KUY}), which represents the interface of a combustion front spreading over time as a closed curve:
\begin{equation}\label{eq:ks_u}
  \dfrac{\partial u}{\partial t} + \dfrac{\delta}{R^4}\dfrac{\partial^4 u}{\partial \sigma^4} + \dfrac{1}{R^2}\left(\alpha - 1 + \dfrac{\delta}{R^2}\right)\dfrac{\partial^2 u}{\partial \sigma^2} + \dfrac{\alpha - 1}{R^2}u - \dfrac{v_c}{2R^2}\left(\dfrac{\partial u}{\partial \sigma}\right)^2 = 0,\\
\end{equation}
where the solution $u: [0, 2 \pi] \times (0, T); (\sigma, t)\mapsto \mathbb{R}$ denotes a height function from an expanding circle solution of \eqref{eq:interfacial} (that is, the solution of \eqref{eq:R} defined below) with radius $R = R(t)$ at time $t$, $\delta$ is a positive parameter, $\alpha > 0$ is a scaled Lewis number, and $v_c$ corresponds to a constant velocity for a uniformly traveling wave described through \eqref{eq:original_ks} below.

\eqref{eq:ks_u} is a perturbed equation of an expanding circle solution of the following time evolution equation:
\begin{equation}\label{eq:interfacial}
  \dfrac{\partial \bs{X}}{\partial t} = V \bs{N} + W \bs{T}, \quad V = v_c + (\alpha - 1)\kappa + \delta \dfrac{\partial^2\kappa}{\partial s^2},
\end{equation}
which is a flow of a family of smooth Jordan curves $\{\mathrm{\Gamma}(t)\}_{0 \le t \le T}$ in the plane $\mathbb{R}^2$.
The solution curve $\mathrm{\Gamma}(t)$ is parameterized by a smooth mapping $\bs{X} = \bs{X}(\sigma, t):\, [0, 2 \pi] \times [0, T] \to \mathbb{R}^2$ such as $\mathrm{\Gamma}(t) = \{\bs{X}(\sigma, t);\, \sigma \in [0, 2 \pi]\}$ (in which the positive direction of $\Gamma(t)$ is counterclockwise). In the second equation in \eqref{eq:interfacial}, $\kappa$ denotes the curvature of $\mathrm{\Gamma}(t)$, and $\partial^2\kappa/\partial s^2$ is the second derivative of $\kappa$ with respect to the arc-length parameter $s=\int_0^\sigma g(\sigma, t) d\sigma$, where $g = |\bs{X}_\sigma|$ is the local length, 
$\bs{X}_\sigma=\partial\bs{X}/\partial\sigma$, and $\partial\kappa/\partial s=g^{-1}\partial\kappa/\partial\sigma$. The velocity of the curve is $\partial \bs{X}/\partial t$, which can be decomposed in the outward normal direction $\bs{N} = - \bs{T}^{\bot}$ and the tangential direction $\bs{T}=\bs{X}_\sigma/|\bs{X}_\sigma|$. 
It is well known that the shape of the curves $\mathrm{\Gamma}(t)$ is determined by the normal velocity $V$ only, and that the tangential velocity $W$ does not affect its shape (see Epstein and Gage~\cite{EG}). The details of the derivation from \eqref{eq:interfacial} to \eqref{eq:ks_u} are provided in subsection~\ref{sec:derivation}.
At a certain scale, \eqref{eq:interfacial} is equivalent to the KS equation:
\begin{equation}\label{eq:original_ks}
  \dfrac{\p f}{\p t} + \delta \dfrac{\p^4 f}{\p x^4} + (\alpha - 1) \dfrac{\p^2 f}{\p x^2} + \dfrac{v_c}{2} \left(\dfrac{\p f}{\p x}\right)^2= 0,
\end{equation}
where $f : [0, L] \times \mathbb{R}_{+}; (x, t) \mapsto \mathbb{R}$ represents a graph of perturbation for a unidirectional uniformly traveling wave solution (see \cite{FS} and \cite{GKY} for details).
From the point of view of the gaseous combustion theory, $\delta$ is considered to be $4$ when $\alpha$ is close to $1$.
We remark that by changing the variables such that $f(x, t) = -u(\sigma, t)$, $x = R\sigma$ with a constant $R$, and taking the limit $R \to \infty$, \eqref{eq:original_ks} can be formally obtained from \eqref{eq:ks_u}. 
In this sense, \eqref{eq:ks_u} corresponds to KS equation \eqref{eq:original_ks}. 

In this paper, to analyze \eqref{eq:ks_u}, we rewrite \eqref{eq:ks_u} using $v = \partial u/\partial \sigma$ as
\begin{equation}\label{eq:ks_v}
  \dfrac{\partial v}{\partial t} + \dfrac{\delta}{R^4}\dfrac{\partial^4 v}{\partial \sigma^4} + \dfrac{1}{R^2}\left(\alpha - 1 + \dfrac{\delta}{R^2}\right)\dfrac{\partial^2 v}{\partial \sigma^2} + \dfrac{\alpha - 1}{R^2} v - \dfrac{v_c}{R^2}v \dfrac{\partial v}{\partial \sigma} = 0,
\end{equation}
and we impose \eqref{eq:ks_v} under the zero average condition $\int_0^{2\pi}v(\sigma, t)d\sigma = 0$ and the periodic boundary condition $v(\sigma, t) = v(\sigma + 2 \pi, t)$.

The KS equation \eqref{eq:original_ks}, which was originally derived as a model of the flame front propagation by Sivashinsky \cite{S}, and independently as the phase turbulence in the reaction-diffusion system by Kuramoto and Tsuzuki \cite{KT}, has been extensively studied for nearly 40 years, both based on the theory of combustion phenomena and in applied mathematics.
From a mathematics perspective, it is well known that \eqref{eq:original_ks} has a rich solution structure including the so-called spatio-temporal chaos (\cite{MS,SM}), and therefore, many researchers of dynamical systems theory, partial differential equation theory, and numerical analysis have been attracted. 
In particular, the existence and uniqueness of the solutions (\cite{AM,NS,T}), the existence and estimates of the Hausdorff and fractal dimensions estimates for an inertial manifold (\cite{NST,R}), the bifurcation structures through a Fourier mode interaction (\cite{AGH1,AGH2,PK}) and numerical studies on for the dynamical behavior (\cite{HN,PS}) are still actively being researched.

In the context of a numerical analysis, in \cite{A1}, Akrivis applied a finite difference scheme to \eqref{eq:original_ks} under a periodic boundary condition.
The method used to prove our result upon convergence to \eqref{eq:ks_v} is based on the idea of \cite{A1}, and is extended to our scheme for \eqref{eq:ks_u}. Akrivis also reported in \cite{A2} a consistent numerical approach to solving \eqref{eq:original_ks} by using a finite element Galerkin method with an extrapolated Crank--Nicolson scheme.
In both cases, a rigorous error analysis was carried out in order to derive the refined error bound.
In addition, numerous other methods have been proposed to find the numerical solutions to \eqref{eq:original_ks} (see \cite{BC} and the references therein).

We emphasize that \eqref{eq:ks_u} and \eqref{eq:ks_v} are moving boundary problems (by contrast, \eqref{eq:original_ks} is formulated in a fixed region), and it is therefore difficult to apply standard dynamical systems theory and analyze the detailed solution structure through a bifurcation analysis.
These facts motivated us to study the behavior of the solution from a numerical aspect.
The results given in Section~\ref{sec:Numerical_wavenumber} of the present paper, guarantee the guidelines for the wavenumber and parameter selection suggested in~\cite{KUY}.

The present paper is organized as follows. 
In section~\ref{sec:derivation}, we show that \eqref{eq:ks_u} is equivalent to an interfacial equation, which was introduced as a combustion model by Frankel and Sivashinsky \cite{FS}.
Our main results are listed in the following section.
In sections~\ref{sec:result_of_CN} and \ref{sec:result_of_N}, convergence of the solutions to the Crank-Nicolson scheme and Newton's method are presented, respectively.
Our main result of the convergence is described in section~\ref{sec:main_result_of_convergence}, which is summarized as
\[
  \max_{0 \le n \le N} \| \bs{u}^n - \bs{U}^n \|_h \le c (k^2 + h^2),
\]
where $\| \cdot \|_h$ denotes the discrete $\mathscr{L}^2$-norm with the space increment $h$, $\bs{u}^n$ is a vectorized solution to \eqref{eq:ks_u} at the $n$-th time step $t^n = nk$ with the time increment $k$, and $\bs{U}^n$ is an approximation solution at the $n$-th step. In section~\ref{sec:proof}, we give the proofs of the main theorems presented in section~\ref{sec:main_result}, that is, show the existence, convergence and uniqueness to our scheme.
In section~\ref{sec:Numerical_algorithm}, the algorithm used by our scheme is given.
In section~\ref{sec:Numerical_experiments}, several numerical experiments of the solution curves are shown.
In the remaining part of section~\ref{sec:Numerical}, we focus on a theoretical linearized stability analysis, particularly the bifurcation theory of the relation between the wavenumber and parameters for moving solution curves and its numerical analysis.
In the final section~\ref{sec:Conclude}, we provide some concluding remarks and areas of future study. 


\section{Main results}\label{sec:main_result}

\subsection{Derivation \eqref{eq:ks_u} from an interfacial equation}\label{sec:derivation}

For convenience, in this subsection, the derivation of \eqref{eq:ks_u} from \eqref{eq:interfacial} for a certain space and time scale is given according to \cite{KUY}.

\eqref{eq:interfacial} has a circle solution
\begin{align}
  \boldsymbol{X}(\sigma, t) &= \boldsymbol{X}_{R} = R(t)\boldsymbol{y}(\sigma), \quad \boldsymbol{y} = \begin{pmatrix} \cos{\sigma}\\ \sin{\sigma}  \end{pmatrix},\label{eq:circle_sol} \\
  \dfrac{dR}{dt} &= v_c + \dfrac{\alpha - 1}{R}.\label{eq:R}
\end{align}
The solution to \eqref{eq:R} is
\begin{equation}\label{eq:R_sol}
  \dfrac{1}{V_{0}}\left\{ R(t) - R(0) - \dfrac{\alpha - 1}{v_c} \log{\dfrac{v_c R(t) + \alpha - 1}{v_c R(0) + \alpha - 1}}\right\} = t,
\end{equation}
which will be used in our numerical scheme (see Section~\ref{sec:Numerical}).
The following property is easily shown.
\begin{proposition}[Proposition~1 in \cite{KUY}]
  Let $\alpha - 1 > 0$.
  Then, for any $R(0) > 0$, the solution $R(t)$ strict monotonically increases, and $R(t) \to +\infty$ holds as $t \to +\infty$.
\end{proposition}

We now consider a perturbation, e.g., $u(\sigma, t)$, of the circle solution \eqref{eq:circle_sol} such that $\boldsymbol{X} - \boldsymbol{X}_R = u(\sigma, t) \boldsymbol{y}$.
By taking the inner product of \eqref{eq:interfacial} and $\boldsymbol{N}$, and denoting $\partial^j \boldsymbol{X}/\partial \sigma^j = F_{j}\boldsymbol{y} + G_{j}{\boldsymbol{y}^{\bot}}$ ($(a, b)^{\bot} = (-b, a)$), \eqref{eq:interfacial} yields
\begin{equation}\label{eq:u}
  \begin{aligned}
    &\dot R + \dot h = \dfrac{v_c g}{G_{1}} + \dfrac{1}{g^{2}}\left(\dfrac{F_{1}G_{2}}{G_{1}} - F_{2}\right)\left\{(\alpha - 1) + \dfrac{3 \delta}{g^{3}}\left(\dfrac{5 g_{\sigma}^{2}}{g} - g_{\sigma\sigma}\right)\right\}\\
    &+ \dfrac{\delta}{g^{4}} \left\{ \dfrac{F_{1}G_{4}}{G_{1}} - F_{4} - \dfrac{7 g_{\sigma}}{g} \left(\dfrac{F_{1}G_{3}}{G_{1}} - F_{3} \right) + \dfrac{F_{2}G_{3} - F_{3}G_{2}}{G_{1}} \right\}.
  \end{aligned}
\end{equation}

We set $\varepsilon = |\alpha - 1| > 0$ as a small parameter and rescale \eqref{eq:u} from $u$ to $f$ using $\varepsilon$ such that $u(\sigma, t) = \varepsilon f(\sigma, \tau)$, $\tau = \varepsilon^2 t$, and $r = \varepsilon^{1/2} R$.
Thus, \eqref{eq:u} can be rewritten as
\begin{align*}
& \Bigg(\dfrac{\partial f}{\partial\tau} + \dfrac{\delta}{r^{4}}
  \dfrac{\partial f^4}{\partial\sigma^4} + \dfrac{1}{r^{2}}\left(\textrm{sgn}(\alpha - 1) + \dfrac{\delta}{r^{2}}\right) \dfrac{\partial f}{\partial\sigma} \\
  &\qquad + \dfrac{\textrm{sgn}(\alpha - 1)}{r^{2}}f - \dfrac{v_c}{2r^{2}}
  \left(\dfrac{\partial f}{\partial\sigma}\right)^{2}\Bigg)\varepsilon^3 + o(\varepsilon^3) = 0.
\end{align*}
Because $\alpha - 1 = \textrm{sgn}(\alpha - 1)\varepsilon$, omitting the term $o(\varepsilon^3)$ and replacing $f$ with $u = u(\sigma, t) \in \mathbb{R}$, one can extract \eqref{eq:ks_u} at the original scale.
Throughout this paper, we assume $u(\sigma, t) \in C^{\infty}$ and $2\pi$-periodic function on $\sigma$, i.e., $u(\sigma, t) = u(\sigma + 2\pi, t)$.
In addition, we set $u(\cdot, 0) = u_0(\sigma)$ and $R(0) = R_0 > 0$ as the initial conditions.
\begin{figure}
  \centering
  \includegraphics[width=0.5\columnwidth]{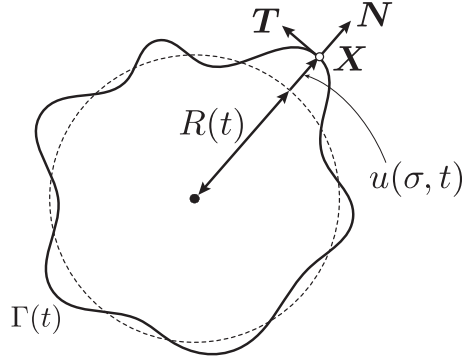}
  \caption{Schematic image of the circle solution to \eqref{eq:interfacial} and a perturbation $u(\sigma, t)$.}
\end{figure}


\subsection{Convergence of solutions to Crank--Nicolson scheme}\label{sec:result_of_CN}
We now state our main results, which will be proved in Section~\ref{sec:proof}.
Let $J \in \mathbb{N}$, $N \in \mathbb{N}$, $T > 0$, $L = 2\pi$, $h = L/J$, $\sigma_i = ih$, $i \in \mathbb{Z}$, $k = T/N$, $t^n = nk$ $(n = 0, 1, \dots, N)$ and
\begin{equation*}
  \mathbb{R}^J_{\mathrm{per}} := \{ \bs{V} = (V_i)_{ i \in \mathbb{Z} };\, V_i \in \mathbb{R} \quad \textrm{and} \quad V_{i+J} = V_i, \quad i\in\mathbb{Z} \}.
\end{equation*}
For $\bs{V} \in \mathbb{R}^J_{\mathrm{per}}$, we describe
\begin{align*}
  \Delta_h V_i &= \dfrac{1}{h^2}(V_{i-1} + 2 V_i + V_{i+1}),\\
  \Delta_h^2 V_i &= \dfrac{1}{h^2}(\Delta_h V_{i-1} - 2 \Delta_h V_i + \Delta_h V_{i+1})
\end{align*}
and for $\bs{V}^0, \dots, \bs{V}^N \in \mathbb{R}^J_{\mathrm{per}}$, we set
\begin{equation*}
  \partial \bs{V}^n = \dfrac{1}{k}(\bs{V}^{n+1} - \bs{V}^n), \quad \bs{V}^{n + \frac{1}{2}} = \dfrac{1}{2}(\bs{V}^n + \bs{V}^{n+1}).
\end{equation*}

We discretize the equation \eqref{eq:ks_v} using the Crank--Nicolson type finite difference scheme.
More precisely, we approximate $\bs{v}^n \in \mathbb{R}^{J}_{\mathrm{per}}$ ($v^n_i = v(\sigma_i, t^n)$) using $\bs{V}^n \in \mathbb{R}^J_{\mathrm{per}}$, where $\bs{V}^0 = \bs{v}^0$, and for $n = 0, 1, \dots, N-1$
\begin{equation}\label{eq:discrete1}
  \begin{multlined}
    \partial \bs{V}^n + \mathcal{L}_h^{n+\frac{1}{2}}\bs{V}^{n + \frac{1}{2}} = \dfrac{v_c}{6 h R_{n+\frac{1}{2}}^2}\bs{\varphi}(\bs{V}^{n + \frac{1}{2}}, \bs{V}^{n + \frac{1}{2}}), \quad i = 1, \dots, J.
  \end{multlined}
\end{equation}
Hereafter, we set $R_{n} = R(t^{n})$, 
\begin{align*}
  & \mathcal{L}_{h}^{n}:\, \mathbb{R}^{J}_{\mathrm{per}} \to \mathbb{R}^{J}_{\mathrm{per}}; \\
  & \quad\left(\mathcal{L}_h^{n}\bs{V}\right)_{i} := \dfrac{\delta}{R_{n}^4}\Delta_h^2 V_{i} + \dfrac{1}{R_n^2}\left(\alpha - 1 + \dfrac{\delta}{R_n^2}\right)\Delta_h V_{i} + \dfrac{\alpha - 1}{R_{n}^2} V_{i},\\
  \intertext{and}
  & \bs{\varphi}:\, \mathbb{R}^J_{\mathrm{per}} \times \mathbb{R}^J_{\mathrm{per}} \to \mathbb{R}^J_{\mathrm{per}}; \\
  & \quad \left(\bs{\varphi}(\bs{V}, \bs{W})\right)_i = (V_{i-1} + V_i + V_{i+1})(W_{i+1} - W_{i-1}). 
\end{align*}
Remark that $R(t)$ is a given function.

The discrete $\mathscr{L}^2$-norm $\| \cdot \|_h$ is introduced in $\mathbb{R}^J_{\mathrm{per}}$ the discrete $\mathscr{L}^2$-norm $\| \cdot \|_h$ using
\begin{equation*}
  \| \bs{V} \|_h := \left( h \sum_{i = 1}^J (V_i)^2 \right)^{\frac{1}{2}}, \quad \bs{V} \in \mathbb{R}_{\mathrm{per}}^J,
\end{equation*}
which is induced through the $\mathscr{L}^2$-inner product $(\cdot, \cdot)_h$ in $\mathbb{R}_{\mathrm{per}}^J$ such that
\begin{equation*}
  (\bs{V}, \bs{W})_h := h \sum_{i = 1}^J V_i W_i, \quad \bs{V}, \bs{W} \in \mathbb{R}_{\mathrm{per}}^J.
\end{equation*}

We can then obtain the result regarding the existence of the numerical solution to \eqref{eq:discrete1}:
\begin{proposition}\label{prop:existence}
  Let $R(0) > \sqrt{\delta/(\alpha - 1)}$ and
  \begin{equation}\label{eq:k}
    k < \dfrac{8\delta}{\left(\alpha - 1 -\dfrac{\delta}{R(T)^2}\right)^2}.
  \end{equation}
  Then, an approximate solution to \eqref{eq:discrete1} exists.
\end{proposition}
As one of the main results in this paper, the following insists that the second-order error estimate to the Crank--Nicolson scheme holds:
\begin{theorem}\label{thm:1}
Let $v$ be sufficiently smooth, and let $M$ be a positive constant such that $\displaystyle\sup_{0 \le t \le T}\left\|\dfrac{\partial v}{\partial \sigma} \right\|_{\infty} \le M$. Suppose that $\bs{V}_{\mathrm{CN}}^1, \dots, \bs{V}_{\mathrm{CN}}^N \in \mathbb{R}_{\mathrm{per}}^J$ are solutions to \eqref{eq:discrete1} with the initial data $\bs{V}_{\mathrm{CN}}^0 = \bs{v}^0$.
  Then, for a sufficiently small $k$, there exists a constant $c = c(\delta, \alpha, R_0, M)$ independent of $k$ and $h$ such that
  \begin{equation}\label{eq:thm1}
    \max_{0 \le n \le N} \|\bs{v}^n - \bs{V}_{\mathrm{CN}}^n \|_h \le c (k^2 + h^2).
  \end{equation}
\end{theorem}
Hereafter, $c$ and $C$ denote general constants that are independent of $k$ and $h$, and do not need to be the same in any two places unless subscripts are applied.
Furthermore, for a sufficiently small $k h^{-\frac{1}{5}}$, the Crank--Nicolson approximations are uniquely defined by \eqref{eq:discrete1}:
\begin{proposition}\label{prop:uniqueness}
  Under the same assumption as in Theorem~\ref{thm:1}, for a sufficiently small $k h^{-\frac{1}{5}}$, the solution to the Crank--Nicolson scheme is unique.
\end{proposition}


\subsection{Convergence of solutions to Newton's method}\label{sec:result_of_N}

To compute the Crank--Nicolson approximations $\bs{V}_{\mathrm{CN}}^1, \dots, \bs{V}_{\mathrm{CN}}^N$, it is necessary to solve a $J \times J$ nonlinear system at each time step.
In this section, we discuss the approximate solution to \eqref{eq:discrete1} using the following Newton's method.

For $n \ge 0$, the linearized system using Newton's method is
\begin{equation}\label{eq:discrete2}
  \begin{aligned}
   \partial \bs{W}^n + \mathcal{L}_h^{n + \frac{1}{2}} \bs{W}^{n + \frac{1}{2}}
    &= \dfrac{v_c}{24 h R_{n + \frac{1}{2}}^2}\bigg(\bs{\psi}(\bs{W}^{n} + \hat{\bs{W}}^{n+1}, \bs{W}^{n+1} - \hat{\bs{W}}^{n + 1})\\
    &+ \bs{\varphi}(\bs{W}^{n} + \hat{\bs{W}}^{n+1}, \bs{W}^{n} + \hat{\bs{W}}^{n + 1})\bigg).
  \end{aligned}
\end{equation}
Here, $\hat{\bs{W}}^0 := \bs{v}^0$, $\hat{\bs{W}}^{n + 1} := 2\bs{W}^n - \bs{W}^{n - 1}$ $(n \ge 1)$ and $\hat{\bs{W}}^1$ is given by
\begin{equation}\label{eq:initial}
  \begin{multlined}
    \partial \hat{\bs{W}}^0 + \mathcal{L}_h^{\frac{1}{2}} \hat{\bs{W}}^{\frac{1}{2}} - \dfrac{v_c}{6 h R_{\frac{1}{2}}^2}\bs{\varphi}(\bs{v}^{0}, \bs{v}^{0})  = 0.
  \end{multlined}
\end{equation}
Note that the symbol `` $\hat{}$ '' indicates the initial value when applying Newton's method.
Although it is still complicated to compute $\bs{W}^{n+1}$, and because \eqref{eq:discrete2} is a $J \times J$ linear system whose matrix depends explicitly on $\bs{W}^n$ and $\bs{W}^{n-1}$, to simplify the problem, we will approximate the solution vectors $\bs{W}^{n}$ through the following process. 

For every time steps, we set $j_n \in \mathbb{N}$ $(n = 1, 2, \dots, N)$ as the maximum number of iterations of Newton's method.
We define the sequences of the approximation vectors as $\bs{W}^{n, j} \in \mathbb{R}_{\mathrm{per}}^J$ $(j = 0, 1, \dots, j_n)$ and $\VN^n := \bs{W}^{n, j_n}$, which corresponds to the approximation solution vector of $\bs{V}_{\mathrm{CN}}^n$.
More precisely, letting $\VN^0 := \bs{V}_{\mathrm{CN}}^0 = \bs{v}^0$, $\hat{\bs{W}}^1 = \hat{\bs{V}}^1 := \VN^1$, $\bs{W}^{n+1, 0} := \hatVN^{n+1}$ and replacing $\bs{W}^{n + 1}$ on the left hand side and $\bs{\psi}$ of \eqref{eq:discrete2} with $\bs{W}^{n+1, j+1}$ and $\bs{W}^{n+1, j}$, respectively, we solve the following:
\begin{equation}\label{eq:discrete3}
  \begin{aligned}
    &\dfrac{1}{k}(\bs{W}^{n+1, j+1} - \VN^n) + \mathcal{L}_h^{n + \frac{1}{2}}(\bs{W}^{n+1, j+1} + \VN^n)\\
    &= \dfrac{v_c}{24 h R_{n + \frac{1}{2}}^2}\big(\bs{\psi}(\VN^n + \hatVN^{n+1}, \bs{W}^{n+1, j} - \hatVN^{n + 1})\\
    &+ \bs{\varphi}(\VN^{n} + \hatVN^{n+1}, \VN^{n} + \hatVN^{n + 1}) \big).
   \end{aligned}
\end{equation}
where $j = 0, \dots, j_{n+1} - 1$.
Then, for $R(0) > \sqrt{\delta/(\alpha - 1)}$ and \eqref{eq:k}, the coefficient matrix of the linear systems \eqref{eq:initial} and \eqref{eq:discrete3} is positive definite, as shown in Proposition~\ref{prop:existence}.

The following result guarantees that the estimate of the form \eqref{eq:thm1} holds, that is, the second-order error estimates for the scheme \eqref{eq:discrete3} are given.
\begin{theorem}\label{thm:2}
Suppose that $\VN^1, \dots, \VN^N \in \mathbb{R}_{\mathrm{per}}^J$ are solutions to \eqref{eq:initial} and \eqref{eq:discrete3} with the initial data $\VN^0 = \bs{v}^0$.
  Under the same assumption as in Theorem~\ref{thm:1}, for a sufficiently small $k$ and $h$ satisfying $k = o(h^{\frac{1}{4}})$, we have
  \begin{equation}
    \max_{0 \le n \le N} \| \bs{v}^n - \VN^n \|_h \le c(k^2 + h^2).
  \end{equation}
\end{theorem}


\subsection{Main results of convergence}\label{sec:main_result_of_convergence}
Let $\bs{V}^n$ be $\bs{V}_{\mathrm{CN}}^n$ or $\VN^n$.
The numerical solution to the integral form \eqref{eq:ks_u}, e.g., $U_i^n$, can be easily calculated from $V_i^n$ as follows. 

We describe $u_i^n = u(\sigma_i, t^n)$ as
\begin{equation}\label{eq:h}
  u_i^n = \dfrac{1}{2 \pi} \int_0^{2\pi} u(\sigma,t^n) d\sigma - \dfrac{1}{2 \pi} \int_0^{2\pi} \left(\int_0^\sigma v(\xi, t^n) d\xi \right) d\sigma + \int_0^{\sigma_i} v(\xi, t^n) d\xi.
\end{equation}
The mean value $I(t) := \frac{1}{2 \pi}\int_0^{2 \pi} u(\sigma, t) d\sigma$ satisfies
\begin{equation}\label{eq:I}
  \dot I(t) = - \dfrac{\alpha - 1}{R(t)^2}I(t) + \dfrac{v_c}{4 \pi R(t)^2}\int_0^{2\pi}v(\sigma, t)^2 d\sigma,
\end{equation}
and therefore
\begin{equation}
  I(t) = \dfrac{\dot R(t)}{\dot R(0)}\left[ I(0) + \dfrac{v_c}{4 \pi \dot R (0)}\int_0^t \left( \dfrac{\dot R (\tau)}{R(\tau)^2} \int_0^{2 \pi} v(\sigma, \tau)^2 d\sigma \right) d\tau \right].
\end{equation}
Here, the symbol `` $\dot{}$ '' indicates $\dot{\mathsf F}(t) = d {\mathsf F}(t)/d t$.

We approximate $u_i^n$ with $U_i^n$ defined by
\begin{equation}\label{eq:U}
  U_i^n := \widetilde I(t^n) - \dfrac{1}{2 \pi} \int_0^{2\pi} \left(\int_0^\sigma \widetilde V(\xi, t^n) d\xi \right) d\sigma + \int_0^{\sigma_i} \widetilde V(\xi, t^n) d\xi,
\end{equation}
where $\widetilde V(\sigma, t)$ is a superposition of the linear interpolation $\widetilde V_i^n(\sigma, t)$ of four numerical values $V_i^n, V_i^{n+1}, V_{i+1}^n, V_{i+1}^{n+1}$, that is,
\begin{equation}\label{eq:interporation}
  \begin{aligned}
    \widetilde V(\sigma, t) &= \sum_{n=0}^{N-1}\sum_{i=0}^{J-1}\widetilde V_i^n(\sigma, t),\\
    \widetilde V_i^n(\sigma, t) &:= (1 - \theta_i) ( (1 - \eta^n)V_{i}^n + \eta^n V_{i}^{n + 1} ) + \theta_i ( (1 - \eta^n) V_{i+1}^{n} + \eta^n V_{i+1}^{n + 1} ),\\ 
    \theta_i &:= \begin{cases}
      (\sigma - \sigma_{i}) h & (\sigma \in [\sigma_i, \sigma_{i+1}])\\ 
      0 & (\mathrm{otherwise})
    \end{cases},\\
    \eta^n &:= \begin{cases}
      (t - t^n) k & (t \in [t^n, t^{n+1}])\\ 
      0 & (\mathrm{otherwise})
    \end{cases},
  \end{aligned}
\end{equation}
and
\begin{equation}\label{eq:hatI}
  \widetilde I(t) = \dfrac{\dot R(t)}{\dot R(0)}\left[ I(0) + \dfrac{v_c}{4 \pi \dot R (0)}\int_0^t \left( \dfrac{\dot R (\tau)}{R(\tau)^2} \int_0^{2 \pi} \widetilde V(\sigma, \tau)^2 d\sigma \right) d\tau \right].
\end{equation}
From direct computations, we find
\begin{equation}
  |I(t^n) - \widetilde I(t^n)| \le c (k^2 + h^2)
\end{equation}
and
\begin{equation}
  \left|\int_0^\sigma \left( v(\xi, t^n) - \widetilde V(\xi, t^n) \right) d\xi \right| \le c (k^2 + h^2).
\end{equation}

Thus, we obtain the main assertion.

\begin{theorem}\label{thm:3}
Let $u$ be sufficiently smooth and $\bs{u}^n \in \mathbb{R}_{\mathrm{per}}^J$ $(u_i^n = u(\sigma_i, t^n))$.
Suppose that $\bs{U}^1, \dots, \bs{U}^N \in \mathbb{R}_{\mathrm{per}}^J$ are given by \eqref{eq:U} with the initial data $\bs{U}^0 = \bs{u}^0$. 
Assume the same for $v = \partial u/\partial \sigma$ as in Theorem \ref{thm:1}.
Then, for a sufficiently small $k$, we have
  \begin{equation}\label{eq:thm3}
    \max_{0 \le n \le N} \|\bs{u}^n - \bs{U}^n \|_h \le c (k^2 + h^2).
  \end{equation}
\end{theorem}


\section{Proofs: Existence, Convergence, and Uniqueness}\label{sec:proof}

In this section, we show the existence of the approximation solutions $\bs{V}^1, \dots, \bs{V}^N \in \mathbb{R}_{\mathrm{per}}^J$ satisfying \eqref{eq:discrete1}, derive the second-order error estimates, and prove the uniqueness of the Crank--Nicolson scheme for a smooth $v$.


\subsection{Preliminaries}\label{sec:Preliminaries}

In addition to the discrete $\mathscr{L}^2$-norm $\| \cdot \|_h$, we use the discrete $\mathscr{H}^1$ and $\mathscr{H}^2$-seminorms, as denoted by $|\cdot|_{1,h}$ and $| \cdot |_{2,h}$, respectively:
\begin{align*}
  | \bs{V} |_{1,h} := \left[ h \sum_{i=1}^J \left( \dfrac{V_{i} - V_{i-1}}{h} \right)^2 \right]^{\frac{1}{2}}, \  | \bs{V} |_{2,h} := \left[ h \sum_{i=1}^J \left( \Delta_h V_i \right)^2 \right]^{\frac{1}{2}}, \ \bs{V} \in \mathbb{R}_{\mathrm{per}}^J.
\end{align*}
Let
\begin{align*}
  & \bs{\psi}:\, \mathbb{R}_{\mathrm{per}}^J \times \mathbb{R}_{\mathrm{per}}^J \to \mathbb{R}_{\mathrm{per}}^J; \\
  &
  \left( \bs{\psi}(\bs{V}, \bs{W}) \right)_i := -(2V_{i-1} + V_i)W_{i-1} + (V_{i+1} - V_{i-1})W_i + (2V_{i+1} + V_i) W_{i+1}.
\end{align*}
We use the following lemma, but omit the proof:
\begin{lemma}[\cite{A1}, Lemma~2.1]\label{lem:1}
  For $\bs{V}, \bs{W}, \bs{U} \in \mathbb{R}^J_{\mathrm{per}}$, we have
  \begin{align}
    &(\bs\varphi(\bs{V}, \bs{W}), \bs{W})_h = - h \sum_{i=1}^J (V_{i+1} - V_{i-2})W_i W_{i-1}, \label{lemma1:1}\\
    &(\bs\varphi( \bs{V}, \bs{V}), \bs{W})_h = - h \sum_{i=1}^J (V_i^2 + V_i V_{i+1} + V_{i+1}^2) (W_{i+1} - W_i), \label{lemma1:2}\\
    &(\bs\psi( \bs{V}, \bs{W}), \bs{U})_h = - h \sum_{i=1}^J \left[ V_i (W_{i+1} + 2W_i) + V_{i+1}(2W_{i+1} + W_i) \right] (U_{i + 1} - U_i), \label{lemma1:3}\\
    &(\varphi(\bs{V}, \bs{V}), \bs{V})_h = 0, \label{lemma1:4}\\
    &\varphi( \bs{V}, \bs{V}) - \varphi(\bs{W}, \bs{W}) = \psi(\bs{W}, \bs{V} - \bs{W}) + \varphi(\bs{V} - \bs{W}, \bs{V} - \bs{W}), \label{lemma1:5}\\
    &-(\Delta_h \bs{V}, \bs{V})_h = | \bs{V} |_{1,h}^2, \label{lemma1:6}\\
    &(\Delta_h^2 \bs{V}, \bs{V})_h = | \bs{V} |_{2,h}^2, \label{lemma1:7}\\
    &| \bs{V} |_{1,h}^2 \le \| \bs{V} \|_h | \bs{V} |_{2,h}, \label{lemma1:8}\\
    &| \bs{V} |_{1,h}^2 \le \eta | \bs{V} |_{2,h}^2 + \dfrac{1}{4\eta}\| \bs{V} \|_h^2, \label{lemma1:9}
  \end{align}
\end{lemma}

Recall that $v(\sigma, t)$ must satisfy $\int_0^{2\pi}v(\sigma, t)d\sigma = 0$ for $t \ge 0$.
If we set $S^n = \sum_{i = 1}^J h V_i^n$, we can see immediately that
\begin{equation}\label{eq:ave_of_V}
  S^{n + 1} = \left(\dfrac{2 R_{n+\frac{1}{2}}^2 - \alpha + 1}{2 R_{n + \frac{1}{2}}^2 + \alpha - 1} \right) S^n
\end{equation}
holds from \eqref{eq:discrete1}.
Thus, if we take the initial condition to satisfy $S^0 = 0$, then $S^n = 0$ holds for all time steps, corresponding to $\int_0^{2\pi} v(\sigma, t) d\sigma = 0$.
Throughout this paper, we assume that $S^n = 0$ $(0 \le n \le N)$ holds.


\subsection{Proof of Proposition~\ref{prop:existence}}\label{sec:proof_existence}

\begin{proof}
  The proof is based on the induction on $n$.
  Assume that $\bs{V}^0, \dots, \bs{V}^n$ $(n < N)$ exist.
  Let $\bs{g}:\, \mathbb{R}_{\mathrm{per}}^J \to \mathbb{R}_{\mathrm{per}}^J$ be defined by
  \begin{align*}
    \bs{g}(\bs{F}) &:= 2 \bs{F} - 2 \bs{V}^n + k \mathcal{L}_h^{n + \frac{1}{2}} \bs{F} - \dfrac{k v_c}{6 h R_{n+\frac{1}{2}}^2} \bs{\varphi}(\bs{F}, \bs{F}).
  \end{align*}
  Taking the inner product with $\bs{F}$ and using \eqref{lemma1:4}, \eqref{lemma1:6}, \eqref{lemma1:7}, \eqref{lemma1:9}, Schwartz's inequality, and the monotonicity of $R(t)$, we have
   \begin{align*}
    (\bs{g}(\bs{F}), \bs{F})_h \ge 2 \| \bs{F} \|_h \left[ \left\{ 1 - \dfrac{k}{8\delta}\left( (\alpha - 1) - \dfrac{\delta}{R(T)^2}  \right)^2 \right\} \| \bs{F} \|_h - \|\bs{V}^n\|_h \right]
  \end{align*}
  under $R(0) > \sqrt{\delta/(\alpha - 1)}$.
  Hence, for
  \begin{equation*}
    k \le \dfrac{8\delta}{\left(\alpha - 1 -\dfrac{\delta}{R(T)^2}\right)^2}
  \end{equation*}
  and
  \begin{align*}
    \| \bs{F} \|_h = \dfrac{8\delta}{8\delta - k\left( (\alpha - 1) - \dfrac{\delta}{R(T)^2}  \right)^2 }\| \bs{V}^n\|_h + 1,
  \end{align*}
  $(\bs{g}(\bs{F}), \bs{F})_h > 0$ holds.
  This yields the existence of $\bs{F}^* \in \mathbb{R}^J_{\mathrm{per}}$ such that $\bs{g}( \bs{F}^* ) = 0$ based on the Brouwer fixed-point theorem (see Lemma~4 in \cite{B}).
  It follows easily that $\bs{V}^{n+1} := 2 \bs{F}^* - \bs{V}^n$ satisfies \eqref{eq:discrete1}.
  
  \smartqed

\end{proof}


\subsection{Proof of Theorem~\ref{thm:1}}\label{sec:proof_convergence}

\begin{proof}
  Let $\bs{r}^n\in\mathbb{R}^J_{\mathrm{per}}$ be the error in the consistency of the scheme \eqref{eq:discrete1}:
  \begin{equation}\label{eq:r}
    \begin{aligned}
      \bs{r}^n &:= \partial \bs{v}^n + \mathcal{L}_h^{n + \frac{1}{2}} \bs{v}^{n + \frac{1}{2}} - \dfrac{V_0}{6 h R_{n+\frac{1}{2}}^2}\bs{\varphi}(\bs{v}^{n+\frac{1}{2}}, \bs{v}^{n+\frac{1}{2}}).
    \end{aligned}
  \end{equation}
  An easy computation shows that
  \begin{equation}\label{eq:maxr}
    \max_{i,n}| r^n_i | \le \max_{n}c(R_n)(k^2 + h^2) \le \exists c (k^2 + h^2).
  \end{equation}
  Set $\bs{e}^n := \bs{v}^n - \bs{V}_{\mathrm{CN}}^n$ $(n=0,\dots,N)$.
  Then, \eqref{eq:discrete1} and \eqref{eq:r} yield
  \begin{equation}\label{eq:e}
    \begin{multlined}
      \partial \bs{e}^n + \mathcal{L}_h^{n + \frac{1}{2}} \bs{e}^{n + \frac{1}{2}}
      + \dfrac{v_c}{6 h R_{n + \frac{1}{2}}^2}\left( \bs{\varphi}(\bs{V}_{\mathrm{CN}}^{n + \frac{1}{2}}, \bs{V}_{\mathrm{CN}}^{n + \frac{1}{2}}) - \bs{\varphi}(\bs{v}^{n + \frac{1}{2}}, \bs{v}^{n + \frac{1}{2}}) \right) = \bs{r}_n.
    \end{multlined}
  \end{equation}
  Here, the non-linear term can be written as
  \begin{align*}
    \bs{\varphi}(\bs{V}_{\mathrm{CN}}^{n + \frac{1}{2}}, \bs{V}_{\mathrm{CN}}^{n + \frac{1}{2}}) - \bs{\varphi}(\bs{v}^{n + \frac{1}{2}}, \bs{v}^{n + \frac{1}{2}})
    &= \bs{\varphi}(\bs{e}^{n + \frac{1}{2}}, \bs{e}^{n + \frac{1}{2}}) - \bs{\varphi}(\bs{e}^{n + \frac{1}{2}}, \bs{v}^{n + \frac{1}{2}})\\
    &- \bs{\varphi}(\bs{v}^{n + \frac{1}{2}}, \bs{e}^{n + \frac{1}{2}}).
  \end{align*}
  Taking the inner product with $\bs{e}^{n + \frac{1}{2}}$, and using \eqref{lemma1:1}, \eqref{lemma1:4}, \eqref{lemma1:6}, \eqref{lemma1:7}, the boundedness of $v_\sigma$, the Schwarz's inequality, and \eqref{eq:maxr}, we have
  \begin{align*}
    &   \|\bs{e}^{n + 1}\|_h - \|\bs{e}^n\|_h \\
    &\le \dfrac{k}{4} \Bigg\{ \dfrac{1}{2\delta}\left( \alpha - 1 - \dfrac{\delta}{R_{n + \frac{1}{2}}^2} \right)^2 + \dfrac{v_c M}{R_{n + \frac{1}{2}}^2} \Bigg\} (\|\bs{e}^{n + 1} \|_h + \|\bs{e}^{n} \|_h)
    + c\sqrt{L}k(k^2 + h^2)\\
    &\le \dfrac{k}{4}\Bigg\{ \dfrac{1}{2\delta}\left( \alpha - 1 - \dfrac{\delta}{R(T)^2} \right)^2 + \dfrac{v_c M}{R(0)^2} \Bigg\}(\|\bs{e}^{n + 1} \|_h + \|\bs{e}^{n} \|_h)
    + c\sqrt{L}k(k^2 + h^2)\\
    &\le C_0 k (\|\bs{e}^{n + 1} \|_h + \|\bs{e}^{n} \|_h) + c\sqrt{L}k(k^2 + h^2), \\
    & C_0 = \dfrac{\alpha - 1}{4 \delta}\left(\dfrac{\alpha - 1}{2} + v_c M \right).  
  \end{align*}
  Here, we used the assumption $R(0) > \sqrt{\delta/(\alpha - 1)}$ and the monotonicity of $R(t)$. 
  Applying Gronwall's discrete inequality we obtain
  \begin{align*}
    \| \bs{e}^{n}\|_h &\le C_1 \| \bs{e}^{n-1}\|_h + \dfrac{c\sqrt{L}k}{1 - C_0 k}(k^2 + h^2)\\ 
    &\le C_1^n \|\bs{e}^0\|_h + \dfrac{C_1^n - 1}{C_1 - 1} \cdot \dfrac{c \sqrt{L}k}{1 - C_0 k}(k^2 + h^2)\\
    &= \dfrac{c \sqrt{L} k}{2 C_0 k} (C_1^n - 1) (k^2 + h^2)\\
    &= C (C_1^n - 1)(k^2 + h^2),   \qquad (n = 1, \dots, N), \\
    C_1 &= \dfrac{1 + C_0 k}{1 - C_0 k}.  
  \end{align*}
  For a sufficiently small $k$, $C_1^n$ can be estimated from above. 
  Indeed, by taking $l \in \mathbb{N}$ satisfying 
  $N \le l \le (1 - C_0 k)/(2C_0 k)$, the following inequality holds. 
  \begin{align*}
    C_1^n = \left( 1 + \dfrac{2C_0 k}{1 - C_0 k} \right)^n \le \left( 1 + \dfrac{2C_0 k}{1-C_0k} \right)^N \le \left( 1 + \dfrac{1}{l} \right)^l < \mathrm{e}.
  \end{align*}
  The proof is completed.
\end{proof}
\qed

\subsection{Proof of Proposition~\ref{prop:uniqueness}}\label{sec:proof_uniqueness}
\begin{proof}
  Suppose $\bs{W}^0 = \bs{v}^0$, and let $\bs{W}^1, \dots, \bs{W}^N \in \mathbb{R}^J_{\mathrm{per}}$ satisfy
  \begin{equation}\label{eq:V}
    \begin{aligned}
      \partial \bs{W}^n + \mathcal{L}_h^{n + \frac{1}{2}} \bs{W}^{n + \frac{1}{2}} = \dfrac{V_0}{6 h R_{n+\frac{1}{2}}^2}\bs{\varphi}(\bs{W}^{n + \frac{1}{2}}, \bs{W}^{n + \frac{1}{2}}).
    \end{aligned}
  \end{equation}
  Substituting $\bs{E}^n := \bs{W}^n - \bs{V}_{\mathrm{CN}}^n$ $(n = 0, \dots, N)$ into \eqref{eq:V} and using \eqref{eq:discrete1} and \eqref{lemma1:5}, we obtain
  \begin{equation}\label{eq:E}
    \begin{aligned}
      \partial \bs{E}^{n} + \mathcal{L}_h^{n + \frac{1}{2}} \bs{E}^{n + \frac{1}{2}} 
       = \dfrac{v_c}{6 h R_{n+\frac{1}{2}}^2}\left( \bs{\psi}(\bs{V}_{\mathrm{CN}}^{n + \frac{1}{2}}, \bs{E}^{n + \frac{1}{2}}) + \bs{\varphi}(\bs{E}^{n + \frac{1}{2}}, \bs{E}^{n + \frac{1}{2}}) \right).
    \end{aligned}
  \end{equation}
  Note that for a sufficiently small $h^{\frac{3}{2}}$, 
  \begin{equation}\label{ineq:U}
    \max_{i,n} |V_i^n| \le \max_{i,n} |v_i^n| + c h \sqrt{h} + c k^2 h^{-\frac{1}{2}} \le C(1 + k^2 h^{-\frac{1}{2}})
  \end{equation}
  follows from \eqref{eq:thm1}.
  Taking in \eqref{eq:E} the inner product with $\bs{E}^{n+\frac{1}{2}}$ and using \eqref{lemma1:3}, \eqref{lemma1:4}, \eqref{lemma1:6}, \eqref{lemma1:7}, \eqref{ineq:U}, and the Schwarz's inequality, we have
  \begin{align*}
    \dfrac{1}{2k}\left( \| \bs{E}^{n + 1}\|_h^2 - \| \bs{E}^{n}\|_h^2 \right) &\le a_1 | \bs{E}^{n+\frac{1}{2}}|_{1,h}^2 - a_2 | \bs{E}^{n + \frac{1}{2}}|_{2,h}^2 - a_3 \| \bs{E}^{n + \frac{1}{2}} \|_h^2\\
    &+ 2 a_4 \| \bs{E}^{n + \frac{1}{2}}\|_h \cdot |\bs{E}^{n + \frac{1}{2}}|_{1,h},
  \end{align*}
  where
  \begin{equation*}
    \begin{aligned}
      a_1 &= \dfrac{1}{R_{n+\frac{1}{2}}^2}\left( \alpha - 1 + \dfrac{\delta}{R_{n+\frac{1}{2}}^2} \right),\quad
      a_2 = \dfrac{\delta}{R_{n+\frac{1}{2}}^4},\quad
      a_3 = \dfrac{\alpha - 1}{R_{n+\frac{1}{2}}^2},\\
      a_4 &= \dfrac{c v_c}{4R_{n+\frac{1}{2}}^2}(1 + k^2 h^{-\frac{1}{2}}).
    \end{aligned}
  \end{equation*}
  Therefore, by using \eqref{lemma1:8}, \eqref{lemma1:9}, we obtain
  \begin{equation}\label{eq:uniq}
    \begin{aligned}
      \|\bs{E}^{n + 1}\|_h^2 - \|\bs{E}^{n}\|_h^2 \le \dfrac{k}{2}\left\{ a_4^2 + \dfrac{(a_1 + 1)^2}{4a_2} - a_3 \right\} (\| \bs{E}^{n + 1} \|_h + \|\bs{E}^n \|_h)^2.
    \end{aligned}
  \end{equation}
  The above inequality implies that for a sufficiently small $k h^{-\frac{1}{5}}$, uniqueness follows immediately through induction.
\end{proof}
\qed

From \eqref{eq:uniq}, we directly obtain the stability result as follows:
\begin{corollary}
  For a sufficiently small $k = \mathcal{O}(h^{\frac{1}{4}})$,
  \begin{equation}
    \|\bs{E}^{n+1}\|_h \le (1 + ck) \|\bs{E}^n\|_h
  \end{equation}
  holds.
\end{corollary}


\subsection{Proof of Theorem~\ref{thm:2}}\label{sec:proof_of_thm2}
\begin{proof}
  We show
  \begin{equation}\label{eq:ass}
    \| \bs{v}^l - \VN^l \|_h^2 \le c_l (k^2 + h^2)^2, \quad l = 0, \dots, N
  \end{equation}
  through an inductive approach.
  Firstly, we can immediately check that \eqref{eq:ass} holds for $l = 0$, and thus $c_0 = 0$.

  Next, assume that \eqref{eq:ass} is valid up to $l = n$, where $c_0 = 0$, $c_1 = 1$, and $\max_{0 \le n \le N}c_n \le c^*$ with a constant $c^*$ independents of $h$ and $k$.
  In the sequel, we additionally assume that $h$ and $k$ are sufficiently small that
  \begin{equation}\label{eq:condition2}
    c^* h^{-1} (k^2 + h^2)^2 \le 1.
  \end{equation}
  We will show later that \eqref{eq:ass} holds for the case of $l = 1$.

  Let $\bs{e}^{n, j} := \bs{v}^n - \bs{W}^{n, j}$ $(j = 0, \dots, j_n)$, $\bs{e}^n := \bs{v}^n - \VN^n$ $(n = 0, \dots, N)$ and $\bs{s}^n \in \mathbb{R}_{\mathrm{per}}^J$ be the consistency error of scheme \eqref{eq:discrete2} such as
  \begin{align}\label{eq:equation_s}
    \bs{s}^n - \bs{r}^n = \dfrac{V_0}{24 h R_{n+\frac{1}{2}}^2} \bs{\varphi}(\bs{v}^{n+1} - \hat{\bs{v}}^{n+1}, \bs{v}^{n+1} - \hat{\bs{v}}^{n + 1}),
  \end{align}
  which follows from \eqref{eq:r} and \eqref{lemma1:5}.
  We can easily see that
  \begin{equation}\label{eq:condition1}
    \max_{i,n}|s_i^n|^2 \le c' (k^2 + h^2)^2.
  \end{equation}
  Using \eqref{eq:r} and \eqref{eq:equation_s} we have
  \begin{equation}\label{eq:4}
    \begin{aligned}
      &\dfrac{1}{k}(\bs{e}^{n+1, j+1} - \bs{e}^n) + \dfrac{1}{2}\mathcal{L}_h^{n + \frac{1}{2}}(\bs{e}^{n+1, j+1} + \bs{e}^n) - \bs{s}^n\\
      &= \dfrac{v_c}{24 h R_{n+\frac{1}{2}}^2} \bigg(\bs{\psi}(\bs{v}^{n+1} - \hat{\bs{v}}^{n+1}, \bs{e}^{n} + \hat{\bs{e}}^{n+1})\\
      &+ \bs{\psi}(\VN^{n} + \hatVN^{n+1}, \bs{e}^{n+1, j} + \bs{e}^n) + \bs{\varphi}(\bs{e}^n + \hat{\bs{e}}^{n+1}, \bs{e}^{n} + \hat{\bs{e}}^{n+1}) \bigg)
    \end{aligned}
  \end{equation}
  from \eqref{eq:discrete3}.
  Taking the inner product with $\bs{e}^{n+1(j+1)} + \bs{e}^n$ and using \eqref{lemma1:2}, \eqref{lemma1:3}, Schwarz's inequality, and the induction hypothesis, we obtain
  \begin{equation}\label{eq:5}
    \begin{aligned}
      &\dfrac{1}{k}\left(\|\bs{e}^{n+1, j+1}\|_h^2 - \|\bs{e}^n\|_h^2 \right) \le a_{n+\frac{1}{2}} |\bs{e}^{n+1, j+1} + \bs{e}^n|_{1,h}^2\\
      &- b_{n+\frac{1}{2}} |\bs{e}^{n+1, j+1} + \bs{e}^n|_{2,h}^2 - c_{n+\frac{1}{2}} \|\bs{e}^{n+1, j+1} + \bs{e}^n\|_h^2\\
      &+ \|\bs{s}^n\|_h \cdot \|\bs{e}^{n+1, j+1} + \bs{e}^n \|_h + d_{n+\frac{1}{2}} |\bs{e}^{n+1, j+1} + \bs{e}^n|_{1,h}\\
      &\times\bigg( M \| \bs{e}^n + \hat{\bs{e}}^{n+1}\|_h + M \| \bs{e}^{n+1 ,j} + \bs{e}^n \|_h + \dfrac{1}{8}h^{-\frac{1}{2}} \|\bs{e}^n + \hat{\bs{e}}^{n+1}\|_h^2 \bigg),
    \end{aligned}
  \end{equation}
  where we set
  \begin{align*}
    a_{n+\frac{1}{2}} &= \dfrac{1}{2R_{n+\frac{1}{2}}^2}\left(\alpha - 1 + \dfrac{\delta}{R_{n+\frac{1}{2}}^2}\right), \quad b_{n+\frac{1}{2}} = \dfrac{\delta}{2R_{n+\frac{1}{2}}^4}, \quad c_{n+\frac{1}{2}} = \dfrac{\alpha - 1}{2 R_{n+\frac{1}{2}}^2},\\
    d_{n+\frac{1}{2}} &= \dfrac{v_c}{R_{n+\frac{1}{2}}^2}
  \end{align*}
  and $M := 1 + \max_{t,\sigma}|v(\sigma,t)|$.
  By the arithmetic-geometric mean inequality, we have
  \begin{align*}
    &d_{n+\frac{1}{2}} M \|\bs{e}^n + \hat{\bs{e}}^{n+1}\|_h \cdot | \bs{e}^{n+1, j+1} + \bs{e}^n |_{1,h} \le \dfrac{d_{n+\frac{1}{2}}^2 M^2}{2 \beta_{n+\frac{1}{2}}} \|\bs{e}^n + \hat{\bs{e}}^{n+1}\|_h^2\\
    &+ \dfrac{\beta_{n+\frac{1}{2}}}{2} | \bs{e}^{n+1, j+1} + \bs{e}^n |_{1,h}^2,\\
    &d_{n+\frac{1}{2}} M \|\bs{e}^{n+1, j} + \bs{e}^{n}\|_h \cdot | \bs{e}^{n+1, j+1} + \bs{e}^n |_{1,h} \le \dfrac{d_{n+\frac{1}{2}}^2 M^2}{2 \beta_{n+\frac{1}{2}}} \|\bs{e}^{n+1, j} + \bs{e}^{n} \|_h^2\\
    &+ \dfrac{\beta_{n+\frac{1}{2}}}{2} | \bs{e}^{n+1, j+1} + \bs{e}^n |_{1,h}^2,\\
    &\dfrac{d_{n+\frac{1}{2}}}{8}h^{-\frac{1}{2}} \|\bs{e}^n + \hat{\bs{e}}^{n+1}\|_h^2 \cdot |\bs{e}^{n+1, j+1} + \bs{e}^n|_{1,h} \le \dfrac{d_{n+\frac{1}{2}}^2}{128\beta_{n+\frac{1}{2}}} \cdot h^{-1} \| \bs{e}^n + \hat{\bs{e}}^{n+1} \|_h^4\\
    &+ \dfrac{\beta_{n+\frac{1}{2}}}{2} | \bs{e}^{n+1, j+1} + \bs{e}^n |_{1,h}^2,\\
    &\left(a_{n+\frac{1}{2}} + \dfrac{3}{2}\beta_{n+\frac{1}{2}} \right) |\bs{e}^{n+1, j+1} + \bs{e}^n |_{1,h}^2
     \le \dfrac{\gamma_{n+\frac{1}{2}}}{2}\|\bs{e}^{n+1, j+1} + \bs{e}^n \|_h^2\\
  &+  \dfrac{1}{2 \gamma_{n+\frac{1}{2}}} \left(a_{n+\frac{1}{2}} + \dfrac{3}{2}\beta_{n+\frac{1}{2}}\right)^2 |\bs{e}^{n+1, j+1} + \bs{e}^n|_{2,h}^2,\\
    &\| \bs{s}^n \|_h \cdot \|\bs{e}^{n+1, j+1} + \bs{e}^n\|_h \le \dfrac{1}{4c_{n+\frac{1}{2}}} \| \bs{s}^n \|_h^2 + c_{n+\frac{1}{2}} \|\bs{e}^{n+1, j+1} + \bs{e}^n \|_h^2,
  \end{align*}
  where $\beta_{n+\frac{1}{2}} = 2a_{n+\frac{1}{2}}/3$ and $\gamma_{n+\frac{1}{2}} = 2a_{n+\frac{1}{2}}^2/b_{n+\frac{1}{2}}$.
   From the above, \eqref{eq:condition1} and \eqref{eq:condition2}, we obtain
  \begin{equation}\label{eq:6}
    \begin{aligned}
      \left( 1 - \tilde{d}_1 k \right) \| \bs{e}^{n+1, j+1} \|_h^2 &\le \tilde{d}_2 k \| \bs{e}^{n + 1, j} \|_h^2 + ( 1 + \tilde{d}_3 k )\| \bs{e}^n \|_h^2\\
      &+ \tilde{d}_4 k \| \bs{e}^{n - 1} \|_h^2 + \tilde{d}_5 k (k^2 + h^2)^2,
    \end{aligned}
  \end{equation}
  where
  \begin{align*}
    \tilde{d}_1 &= \dfrac{2a_{n+\frac{1}{2}}^2}{b_{n+\frac{1}{2}}}, \quad \tilde{d}_2 = \dfrac{3 d_{n+\frac{1}{2}}^2 M^2}{2a_{n+\frac{1}{2}}}, \quad \tilde{d}_3 = \dfrac{2a_{n+\frac{1}{2}}^2}{b_{n+\frac{1}{2}}} + \dfrac{15d_{n+\frac{1}{2}}^2 M^2}{a_{n+\frac{1}{2}}} + \dfrac{135d_{n+\frac{1}{2}}^2}{32a_{n+\frac{1}{2}}},\\
    \tilde{d}_4 &= \dfrac{3 d_{n+\frac{1}{2}}^2 M^2}{2a_{n+\frac{1}{2}}} + \dfrac{15d_{n+\frac{1}{2}}^2}{32a_{n+\frac{1}{2}}}, \quad \tilde{d}_5 = \dfrac{c'}{4c_{n+\frac{1}{2}}}.
  \end{align*}
  Now, we let $D$ be such that $\|\bs{v}^n - \hat{\bs{v}}^n\|_h^2 \le D k^4$, and $p$ and $\tilde p$ be such that, for a sufficiently small $k$, 
  \begin{equation}
    \dfrac{\tilde{d}_2 k}{1 - \tilde{d}_1 k} \le \tilde p k, \quad \dfrac{\delta_{3j} + \tilde{d}_j k}{(1 - \tilde{d}_1 k) (1 - \tilde p k)} \le \delta_{3j} + pk \quad (j = 3,4,5),
  \end{equation}
  where $\delta_{3j}$ is the Kronecker's delta.
  Then, for a sufficiently small $k$, \eqref{eq:6} is rewritten as
  \begin{align*}
    \| \bs{e}^{n+1, j_{n+1}} \|_h^2 &\le c_{n+1} (k^2 + h^2)^2,\\
  \end{align*}
  where
  \begin{align*}
    c_{n+1} &= \left( 2(\tilde p k)^{j_{n+1}} (Dk^4 + 8 c_n + 2 c_{n-1} ) + (1 + pk)c_n + pkc_{n-1} + pk\right).
  \end{align*}

  Finally, we show that \eqref{eq:ass} for $l = 1$.
  Substituting $n = 0$ into \eqref{eq:5}, we have
  \begin{equation}\label{eq:10}
    \begin{aligned}
      \dfrac{1}{k} \| \bs{e}^1 \|_h^2 &\le a_{\frac{1}{2}} |\bs{e}^1|_{1,h}^2 - b_{\frac{1}{2}} |\bs{e}^1|_{2,h}^2 - c_{\frac{1}{2}} \|\bs{e}^1\|_h^2 + \| \bs{s}^0 \|_h \cdot \| \bs{e}^1 \|_h\\
      &+ d_{\frac{1}{2}} \left( M \| \hat{\bs{e}}^1 \|^2_h + M \| \bs{e}^1 \|_h^2 + \dfrac{1}{8} h^{-\frac{1}{2}} \| \hat{\bs{e}}^1 \|_h^2\right)|\bs{e}^1|_{1,h}^2 \\
      &\le  \left(\dfrac{a_{\frac{1}{2}}^2}{b_{\frac{1}{2}}} + \dfrac{3d_{\frac{1}{2}}^2 M^2}{4a_{\frac{1}{2}}}\right) \| \bs{e}^{1} \|_h^2 + \dfrac{3d_{\frac{1}{2}}^2 M^2}{4a_{\frac{1}{2}}} \| \hat{\bs{e}}^{1}\|_h^2 + \dfrac{3d_{\frac{1}{2}}^2}{256a_{\frac{1}{2}}}h^{-1} \| \hat{\bs{e}}^{1}\|_h^4\\
      &+ \dfrac{1}{4c_{\frac{1}{2}}}\| \bs{s}^0 \|_h^2.
    \end{aligned}
  \end{equation}
  Now, let $\hat{\bs{e}}^1 := \bs{v}^1 - \hatVN^1$.
  Using \eqref{eq:initial} and \eqref{eq:r}, we obtain
  \begin{align}
    \partial\hat{\bs{e}}^1 + \dfrac{k}{2} \mathcal{L}_h^{\frac{1}{2}} \hat{\bs{e}}^1
     = \dfrac{v_c k}{24 h R_{\frac{1}{2}}^2}\left( \bs{\psi}(2\bs{v}^0, \bs{v}^1 - \bs{v}^0) + \bs{\varphi}(\bs{v}^1 - \bs{v}^0, \bs{v}^1 - \bs{v}^0) \right) + k \bs{r}^0.
  \end{align}
  Taking the inner product with $\hat{\bs{e}}^1$, and using \eqref{lemma1:6}, \eqref{lemma1:7}, \eqref{lemma1:9}, and the fact that $\left( \bs{\psi}(\bs{V}, \bs{W})\right)_i = (V_i + 2 V_{i+1}(W_{i+1} - W_{i-1}) + (2W_{i-1} + W_i)(V_{i+1} - V_{i-1})$, we obtain    
  \begin{equation}\label{eq:e1}
    \begin{aligned}
      &\| \hat{\bs{e}}^1 \|_h^2 + \dfrac{\delta k}{2 R_{\frac{1}{2}}^4}| \hat{\bs{e}}^1 |_{2,h}^2 - \dfrac{k}{2R_{\frac{1}{2}}^2}\left( \alpha - 1 + \dfrac{\delta}{R_{\frac{1}{2}}^2} \right) | \hat{\bs{e}}^1 |_{1,h}^2 + \dfrac{(\alpha - 1)k}{2 R_{\frac{1}{2}}^2} \| \hat{\bs{e}}^1 \|_h^2\\
      &\le \dfrac{v_c k}{24 h R_{\frac{1}{2}}^2} \left(\| \bs{\psi}(2 \bs{v}^0, \bs{v}^1 - \bs{v}^0) \|_h + \| \bs{\varphi}(\bs{v}^1 - \bs{v}^0, \bs{v}^1 - \bs{v}^0)\|_h \right)\| \hat{\bs{e}}^1 \|_h\\
      & + k \| \bs{r}^0 \|_h \cdot \| \hat{\bs{e}}^1 \|_h.
    \end{aligned}
  \end{equation}
  Here, we remark that the following hold:
  \begin{equation*}
    \| \bs{\psi}(2\bs{v}^0, \bs{v}^1 - \bs{v}^0) \|_h \le c k h, \quad \| \bs{\varphi}(\bs{v}^1 - \bs{v}^0, \bs{v}^1 - \bs{v}^0) \|_h \le c k^2 h.
  \end{equation*}
  Using the above estimations, from \eqref{eq:e1} we have
  \begin{equation}
    \| \hat{\bs{e}}^1 \|_h^2 \le a_{\frac{1}{2}} k | \hat{\bs{e}}^1 |_{1,h}^2 - b_{\frac{1}{2}} k | \hat{\bs{e}}^1 |_{2,h}^2 - c_{\frac{1}{2}} k \| \hat{\bs{e}}^1 \|_h^2 + \dfrac{c v_c k^2}{R^2}\| \hat{\bs{e}}^1 \|_h + k \| \bs{r}^0 \|_h \cdot \| \hat{\bs{e}}^1 \|_h.
  \end{equation}
  This yields
  \begin{align*}
    \| \hat{\bs{e}}^1 \|_h^2 &\le a_{\frac{1}{2}} k \| \hat{\bs{e}}^1 \|_h \cdot | \hat{\bs{e}}^1 |_{2,h} - b_{\frac{1}{2}} k | \hat{\bs{e}}^1 |_{2,h}^2 - c_{\frac{1}{2}} k \| \hat{\bs{e}}^1 \|_h^2 + \dfrac{c v_c k^2}{R^2} \| \hat{\bs{e}}^1 \|_h\\
    &+ k \left(\dfrac{1}{4c_{\frac{1}{2}}}\| \bs{r}^0 \|_h^2 + c_{\frac{1}{2}} \| \hat{\bs{e}}^1 \|_h^2 \right)\\
     &\le \dfrac{a_{\frac{1}{2}}^2k}{4 b_{\frac{1}{2}}} \| \hat{\bs{e}}^1 \|_h^2 + \dfrac{c v_c^2 k^4}{2 R^4} + \dfrac{1}{2} \| \hat{\bs{e}}^1 \|_h^2 + \dfrac{k}{4c_{\frac{1}{2}}}\| \bs{r}^0 \|_h^2.
   \end{align*}
  That is, we have
  \begin{align*}
    \left(\dfrac{1}{2} - \dfrac{a_{\frac{1}{2}}^2}{4b_{\frac{1}{2}}}k\right) \|\hat{\bs{e}}^1\|_h^2 \le \dfrac{c v_c^2 k^4}{2R^4} + \dfrac{k}{4c_{\frac{1}{2}}}(k^2 + h^2)^2.
  \end{align*}
  This implies that, for a sufficiently small $k$,
  \begin{equation}
    \|\hat{\bs{e}}^1\|_h^2 \le C(k^2 + h^2)^2
  \end{equation}
  holds.
  Therefore, \eqref{eq:10} is transformed into
  \begin{equation}
    \begin{aligned}
      \left\{ 1 - k \left(\dfrac{a_{\frac{1}{2}}^2}{b_{\frac{1}{2}}} + \dfrac{3 d_{\frac{1}{2}}^2 M^2}{4a_{\frac{1}{2}}}\right) \right\} \| \bs{e}^1 \|_h^2
      &\le  \dfrac{k}{4}\left[ \dfrac{3 d_{\frac{1}{2}}^2}{a_{\frac{1}{2}}} \left( M^2 + \dfrac{1}{64} \right) \| \hat{\bs{e}}^1 \|_h^2 + \dfrac{1}{c_{\frac{1}{2}}}\|\bs{s}^0\|_h^2 \right].
    \end{aligned}
  \end{equation}
  From the above, for a sufficiently small $k$, we finally obtain
  \begin{equation}
    \| \bs{e}^1 \|_h^2 \le C k ( \| \hat{\bs{e}}^1 \|_h^2 + \| \bs{s}^0 \|_h^2 ) \le c_1 (k^2 + h^2)^2 \quad (c_1 := 1),
  \end{equation}
  which is our assertion.
  Note that it is obvious that, for a sufficiently small $k$, $\displaystyle\max_{0 \le n \le N} c_n \le c^*$ holds.  
\end{proof}
\qed


\section{Our scheme and numerical experiments}\label{sec:Numerical}

\subsection{Algorithm (proposed scheme)}\label{sec:Numerical_algorithm}

We can summarize our numerical scheme using the following symbols. 
\begin{center}
  \begin{tabular}{c|l}\hline
    $\bs{V}_{\mathrm{CN}}^{n}$       & the Crank--Nicolson approximations\\ 
    $\bs{W}^{n, j}$    & the Newton iteration ($0 \le j \le j_{n}$)\\
    $\hat{\bs{W}}^{n}$ & the initial value of the Newton iteration \eqref{eq:discrete2} \\
    $\hatVN^{n}$    & the initial value of \eqref{eq:discrete3} \\
    $\VN^{n}$          & the approximation solution of $\bs{V}_{\mathrm{CN}}^{n}$ by  \eqref{eq:discrete3} \\ \hline
  \end{tabular}
\end{center}

The algorithm is as follows. 
We compute 
$\VN^{0}, \VN^{1}, \VN^{2}, \dots$ for the approximate solution to the differential form \eqref{eq:ks_v}, and 
$\bs{U}^{0}, \bs{U}^{1}, \bs{U}^{2}, \dots$ for the approximate solution to the integral form \eqref{eq:ks_u}. 
Note that true approximate solution $\bs{V}^n$ is implemented as $\VN^{n}$. 

\begin{description}
\item[\textbf{Step 0. Set the initial values}]\mbox{}
  \begin{description}
  \item[Step 0-1] 
    Set the initial values $R_{0}$ and $U_{i}^{0} := u(\sigma_{i}, 0)$; 

  \item[Step 0-2] 
    Set the initial values for \eqref{eq:initial} such that
    $\VN^{0} := \bs{v}^{0} = \dfrac{\partial\bs{u}^{0}}{\partial\sigma}$;  

  \item[Step 0-3]
    Compute $\widetilde I(0)$ by \eqref{eq:hatI};
  \end{description}

\item[\textbf{Step 1. Compute $\VN^1$ and $\bs{U}^{1}$}]\mbox{}
  \begin{description}
  \item[Step 1-1]
    Compute $\hat{\bs{W}}^{1}$ using \eqref{eq:initial}, and put $\VN^1 := \hat{\bs{W}}^{1}$;

  \item[Step 1-2]
    Compute $\widetilde I(t^1)$ by \eqref{eq:hatI}, and obtain $\bs{U}^{1}$ through \eqref{eq:U};
    Put $n = 2$;
  \end{description}

\item[\textbf{Step 2. Compute $\VN^n$ and $\bs{U}^{n}$}]\mbox{} 
  \begin{description}
  \item[Step 2-1]
    Set $\hatVN^{n} = 2 \VN^{n-1} - \VN^{n-2}$; 

  \item[Step 2-2]
    Iterate \eqref{eq:discrete3} for $0 \le j \le j_{n} - 1$; 

  \item[Step 2-3]
    Set $\VN^{n} := \bs{W}^{n, j_{n}}$; 

  \item[Step 2-4]
    Compute $\widetilde I(t^{n})$ using \eqref{eq:hatI}, and obtain $\bs{U}^{n}$ by \eqref{eq:U}; 
  \end{description}
\item[\textbf{Step 3.}] Put $n := n + 1$. 
  If $n\le N$, go to \textbf{Step 2}, 
  else \textbf{END}. 
\end{description}
\begin{remark}
  In Step 0-2, one can use the difference 
  \[
  V_{\mathrm{Newton},i}^{0} := \dfrac{U_{i+1}^{0} - U_{i-1}^{0}}{2 h}, 
  \]
  instead of $\dfrac{\partial\bs{u}^{0}}{\partial\sigma}$. 
  Note that the convergence result holds, 
  even in this case. 
\end{remark}


\subsection{Numerical experiments}\label{sec:Numerical_experiments}

\begin{figure}[h]
  \centering
  \begin{tabular}{@{}cc@{}}
  \subfloat[][]{\includegraphics[width=0.469\columnwidth]{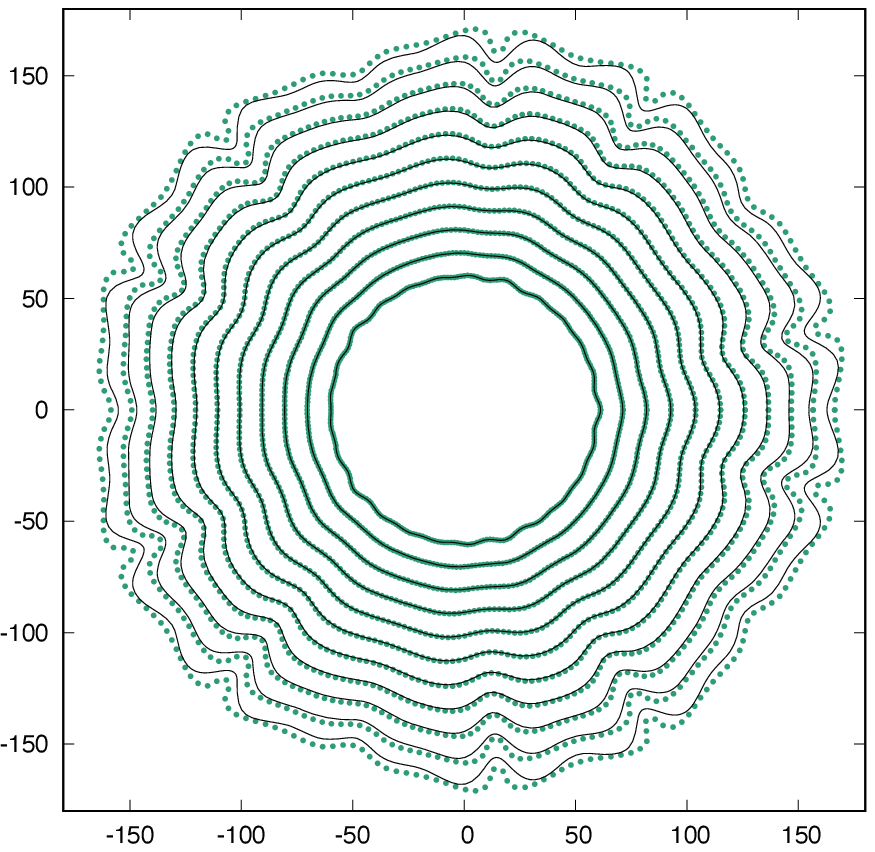}}&
  \subfloat[][]{\includegraphics[width=0.48\columnwidth]{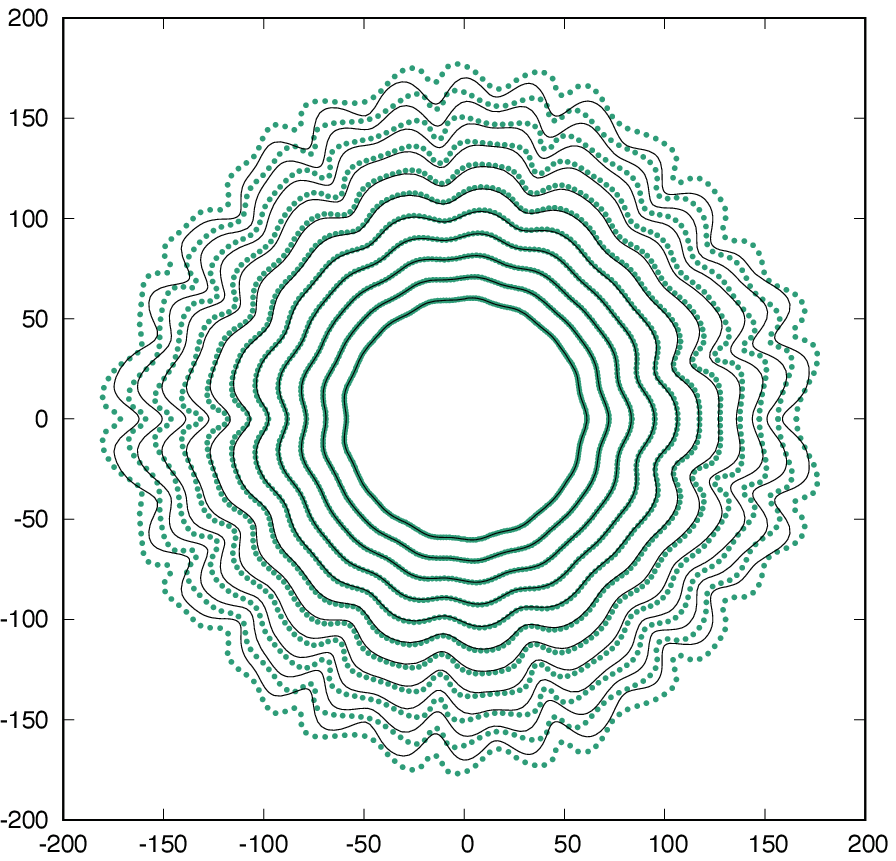}}
  \end{tabular}
  \caption[]{Comparison results for \eqref{eq:ks_u} using our scheme \S\ref{sec:Numerical_algorithm} and for \eqref{eq:interfacial} using the method proposed in \cite{GKY}.
    (a): $p_i = 0.12 \cdot i$, $m_i = 5(1 + i)$;
    (b): $p_i = 0.06 \cdot i$, $m_1 = 7$, $m_2 = 11$, $m_3 = 13$, $m_4 = 17$;
  }\label{fig:1}
\end{figure}

In this subsection, we show the numerical results on \eqref{eq:ks_u} given by the previous subsection and those on the evolution equation for a closed curve \eqref{eq:interfacial} given by the numerical scheme proposed in \cite{GKY}.
In \cite{GKY}, a numerical scheme for \eqref{eq:interfacial} was introduced such that the tangential velocity $W$ controls the grid-point spacing to be uniform, or, more precisely, an asymptotically uniform tangential velocity is chosen according to \cite{BKKMP,SY}.
The scheme was verified based on measurements of the experimental orders of convergence, and even for the initial nontrivial conditions, is of approximately a first order (\cite{GKUY}).
In Fig.~\ref{fig:1}, the green dots and solid lines corresponds to the numerical solution to \eqref{eq:interfacial} and \eqref{eq:ks_u}, respectively.
The solid lines are given by
\begin{equation}\label{eq:circle}
  ( R(t) + u(\sigma, t) )\begin{pmatrix} \cos\sigma \\ \sin\sigma \end{pmatrix}.
\end{equation}
The parameters are $k = 0.01$, $N=1024$, $T = 1000$, $v_c = 0.1$, $\alpha = 1.28$, and $\delta = 4.0$.
The initial conditions are $R(0) = 60.0$ and
\begin{equation}\label{eq:initial_condition}
  u(\sigma, 0) = \sum_{i = 1}^{4}p_{i}\cos{m_i \sigma}.
\end{equation}  
Note that the green dots in Fig.~\ref{fig:1} are drawn at every $4$ points, that is, $4 l$ ($0 \le l \le 256$).
The solution curves are depicted by $t = 0, 100, \dots, 1000$.

From Fig.~\ref{fig:1}, for the perturbation part, it is experimentally confirmed that the behavior of the graph and interface are almost the same.
Therefore, we confirm that it is reasonable to use the graph \eqref{eq:ks_u} to theoretically analyze the qualitative properties of the solution to the interfacial equation \eqref{eq:interfacial}.
There is a difference in the behavior when a sufficient amount of time has been spent, which may indicate the limitation in that the graph cannot be overhang.


\subsection{Wavenumber selection and parameters (theoretical results)}\label{sec:known_results}

\begin{figure}[h]
  \centering
  \includegraphics[width=0.48\columnwidth]{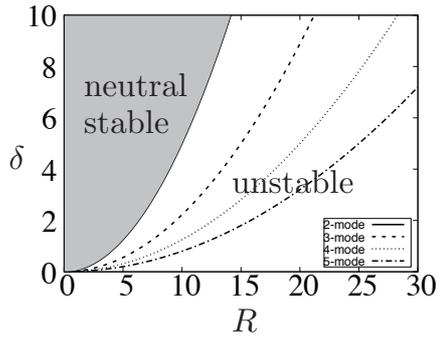}
  \caption{The $|m|$-mode neutral stability curves when $\alpha = 1.2$ with $|m| = 2, 3, 4, 5$.
    The horizontal and vertical axes represent $R$ and $\delta$, respectively (\cite{KUY}).}\label{fig:neutral}
\end{figure}

In this subsection only, let $R(t)$ be the bifurcation parameter $R$.
Then, the neutral stability curves upon which the linearized operator
\begin{align*}
  \mathcal{L} &= \dfrac{\delta}{R^4} \dfrac{\partial^4}{\partial \sigma^4} + \dfrac{1}{R^2}\left(\alpha - 1 + \dfrac{\delta}{R^{2}}\right) \dfrac{\partial^2}{\partial \sigma^2} + \dfrac{\alpha - 1}{R^{2}}
\end{align*}
corresponding to \eqref{eq:ks_u} has eigenvalues of zero around the circle solution are defined as follows (see Fig.~\ref{fig:neutral}):
\begin{definition}[Definition~2 in \cite{KUY}]
  The neutral stability curves are defined as a set of parameters
  \begin{equation}\label{eq:delta-R}
    \left\{(R, \delta);\, \delta = \dfrac{\alpha - 1}{m^2}R^2, \ m \in \mathbb{Z} \right\}
  \end{equation}
  upon which $\mathcal{L}$ has eigenvalues of zero.
\end{definition}
From the above definition, we see that the circle solution is neutrally stable as in the gray region, whereas in the white region where $R > R^* := 2\sqrt{\delta/(\alpha - 1)}$, the circle solution is unstable except at the 2-mode neutral stability curve because at least one eigenvalue is positive (for details, see Appendix~\ref{sec:apA}).
According to this definition, as long as $|\alpha - 1|$ is sufficiently small, we can know the relationship between the wavenumber of the solution to \eqref{eq:interfacial} and the parameters.
Therefore, to study the qualitative properties of the solution for \eqref{eq:interfacial}, it is important to study the instability of the solution for \eqref{eq:ks_u}.


\subsection{Wavenumber selection and parameters (numerical results)}\label{sec:Numerical_wavenumber}
\begin{figure}[h]
  \centering
  \subfloat[][]{\includegraphics[width=0.48\columnwidth]{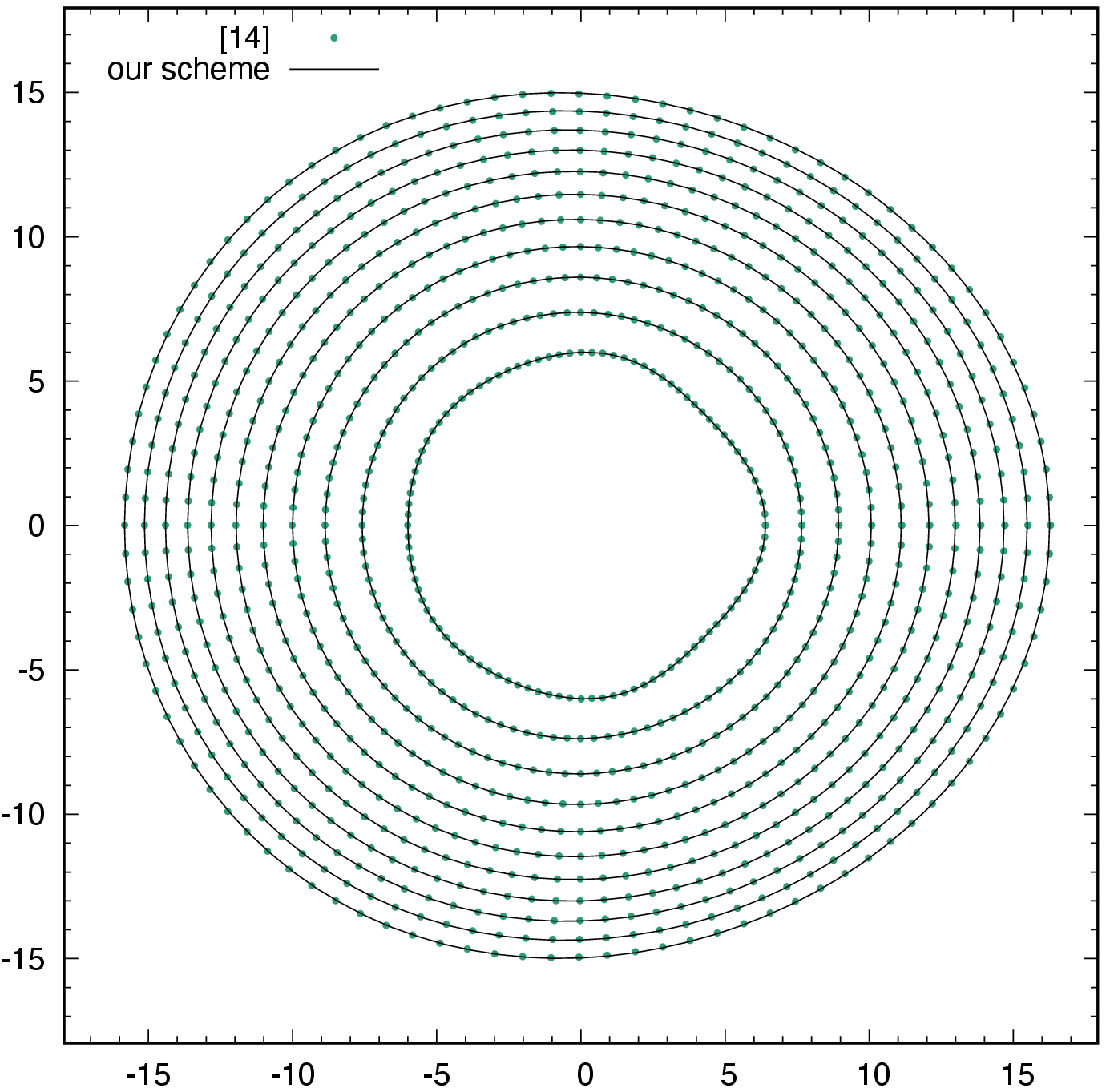}} \\ 
  \subfloat[][]{\includegraphics[width=0.35\columnwidth]{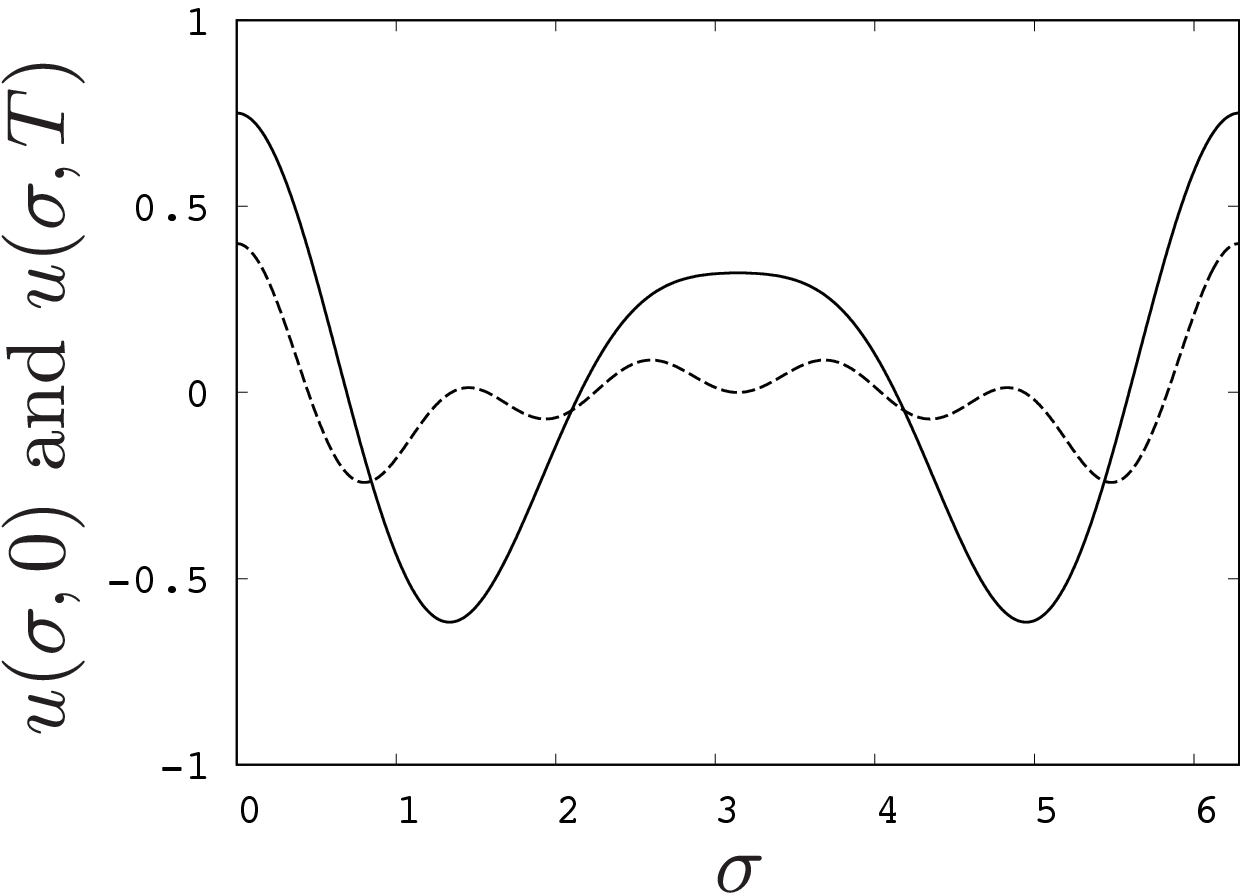}} \hspace{1cm}  
  \subfloat[][]{\includegraphics[width=0.35\columnwidth]{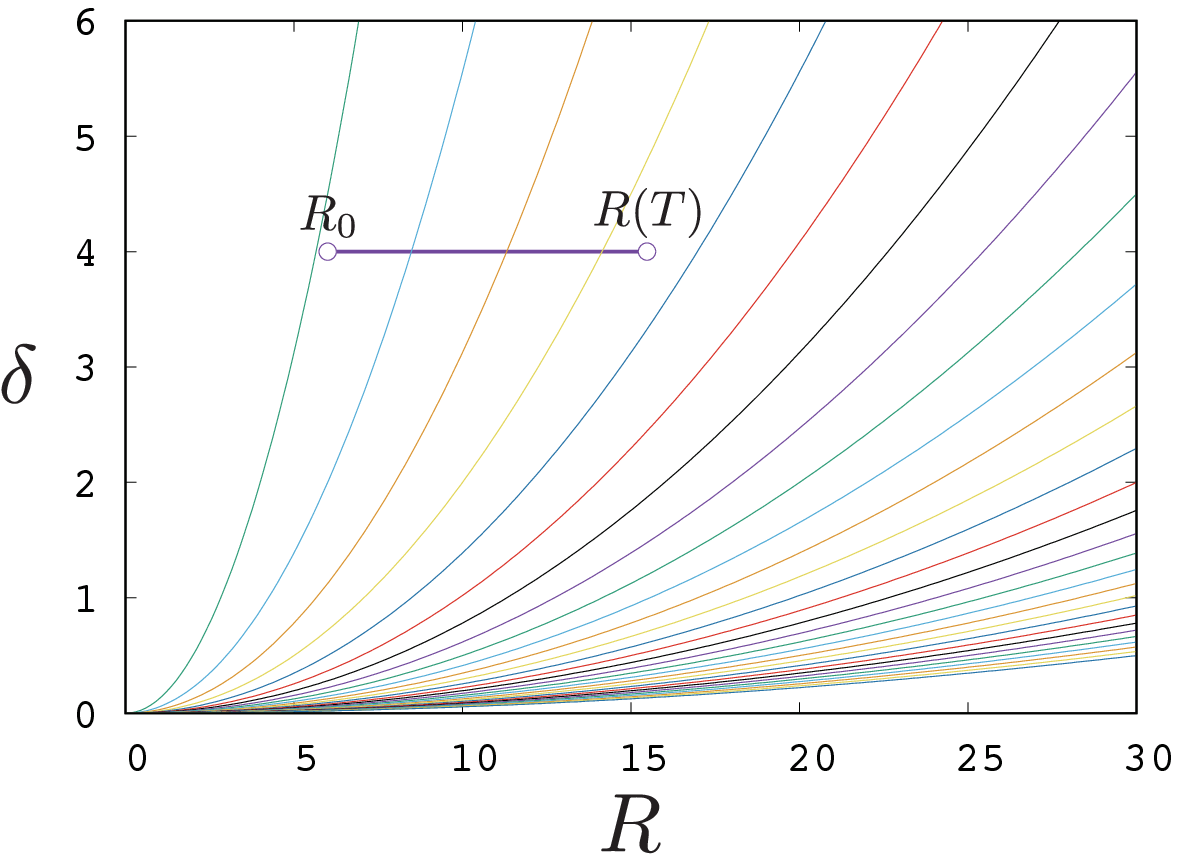}}
  \caption[]{Numerical results for the case in which $R_0 = 6.0$ and $(m_1, m_2, m_3, m_4) = (2, 3, 4, 5)$ in \eqref{eq:initial_condition}.}\label{fig:2}
\end{figure}

\begin{figure}[h]
  \centering
  \subfloat[][]{\includegraphics[width=0.48\columnwidth]{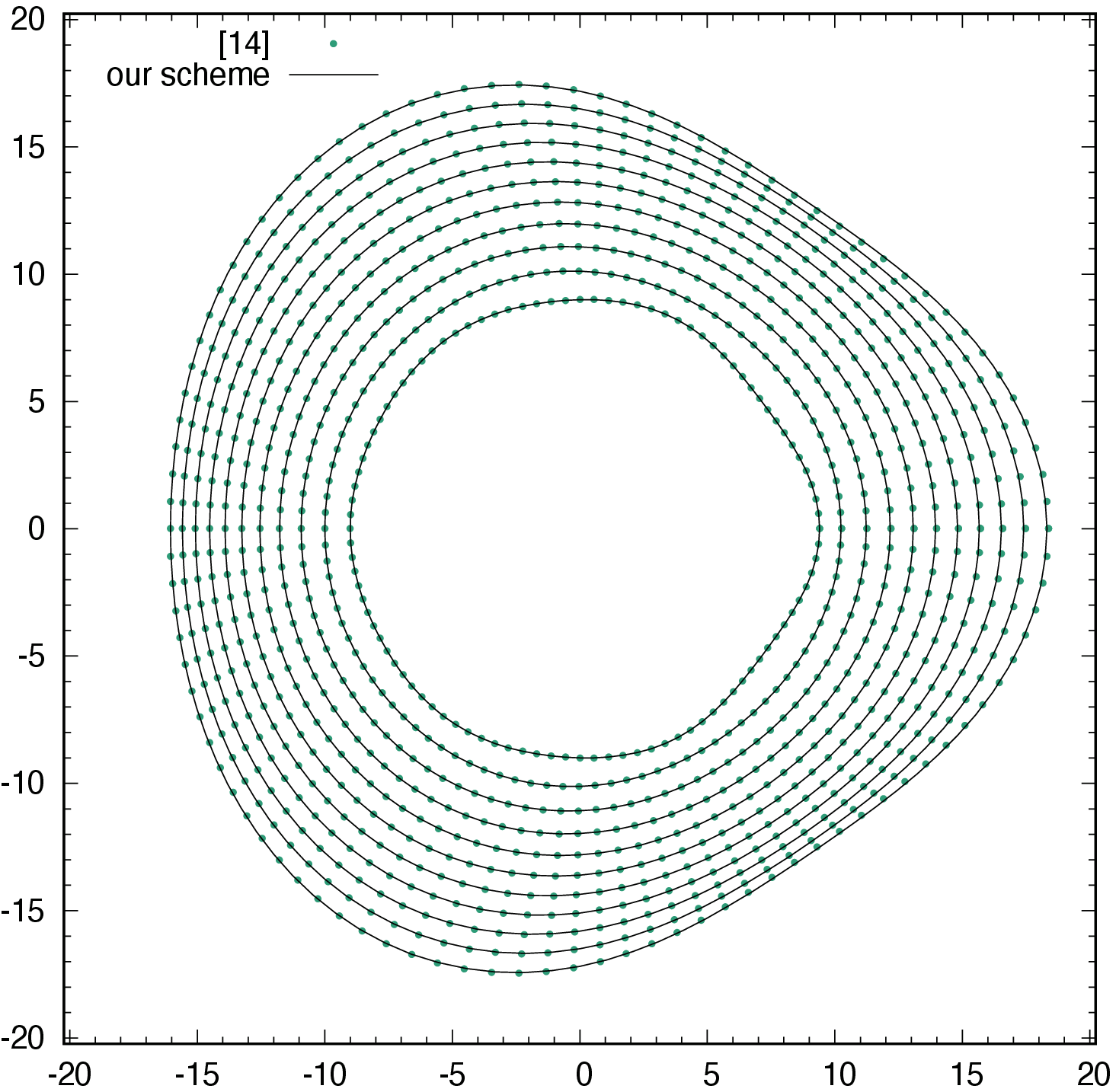}} \\ 
  \subfloat[][]{\includegraphics[width=0.35\columnwidth]{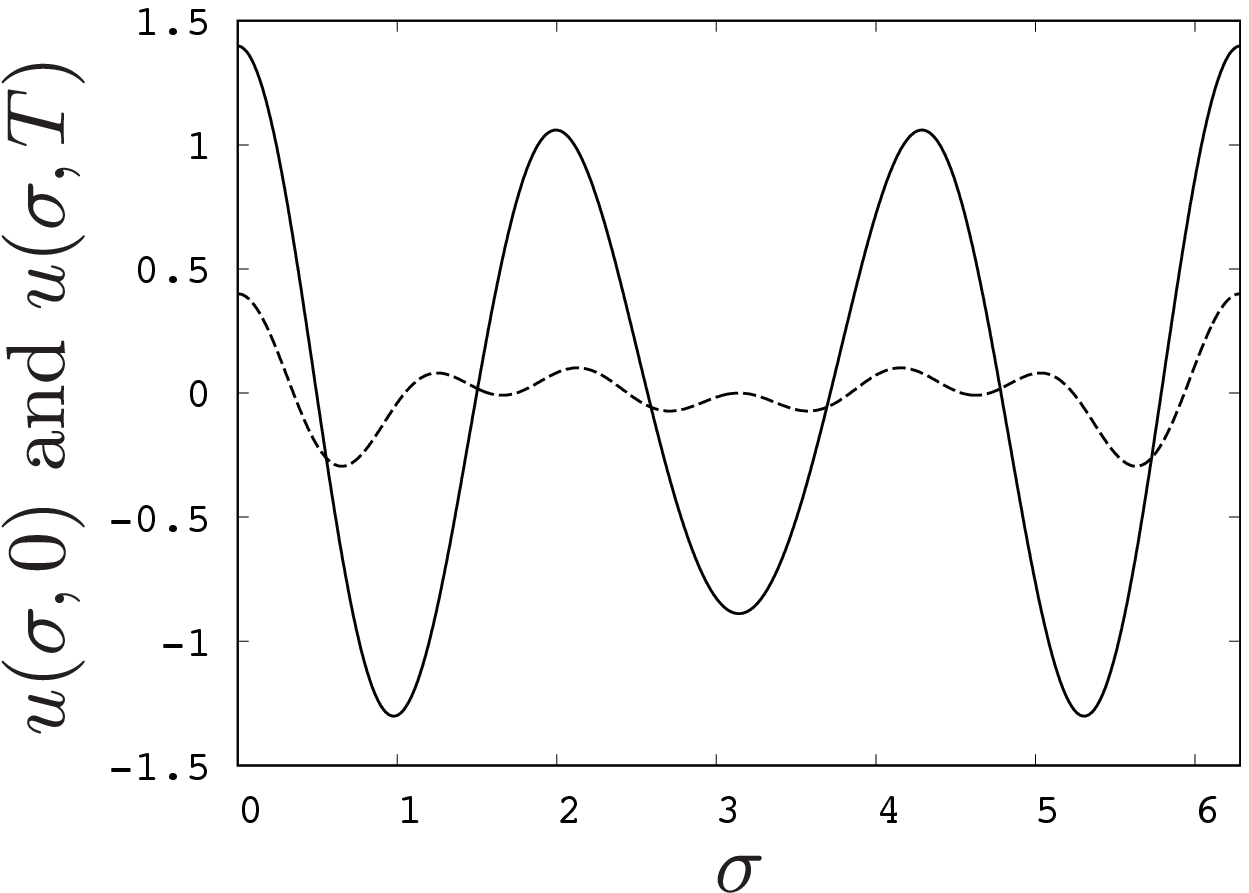}} \hspace{1cm}
  \subfloat[][]{\includegraphics[width=0.35\columnwidth]{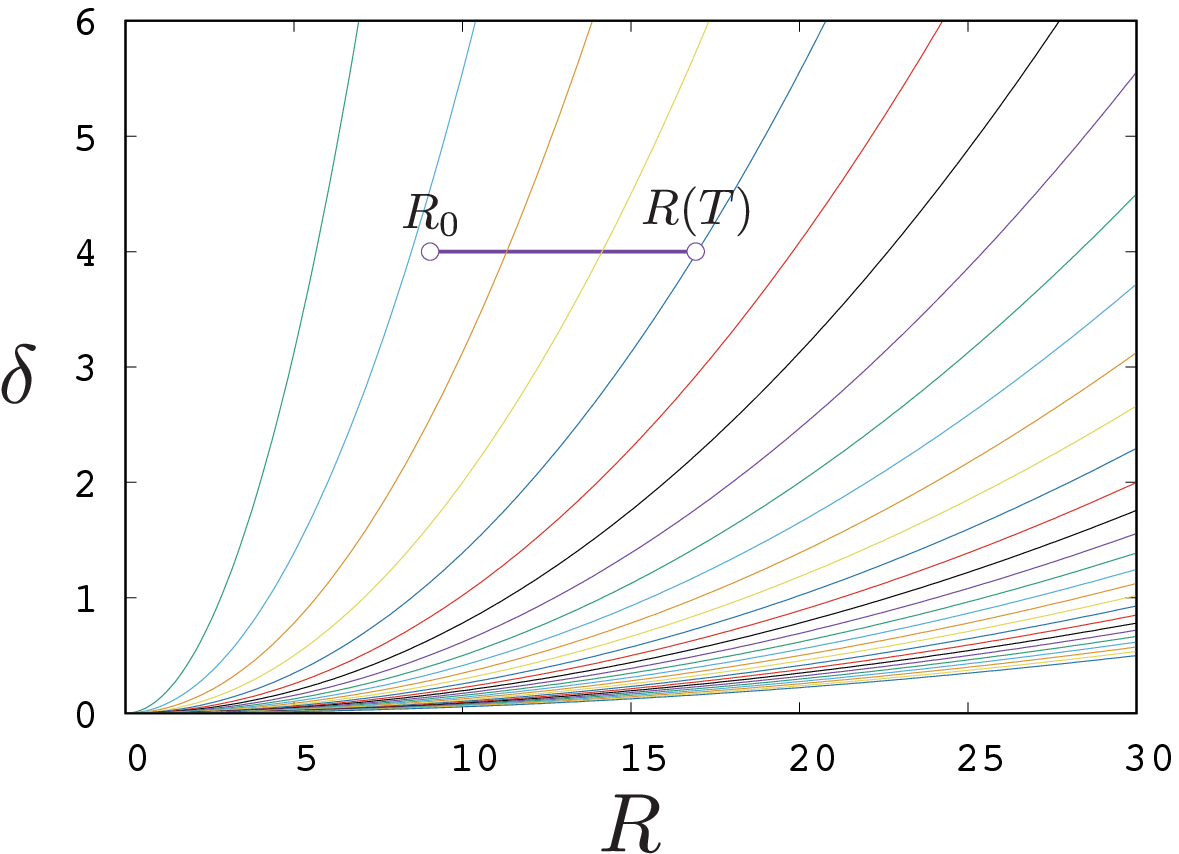}}
  \caption[]{Numerical results for the case in which $R_0 = 9.0$ and $(m_1, m_2, m_3, m_4) = (3, 4, 5, 6)$ in \eqref{eq:initial_condition}.}\label{fig:3}
\end{figure}

\begin{figure}[h]
  \centering
  \subfloat[][]{\includegraphics[width=0.48\columnwidth]{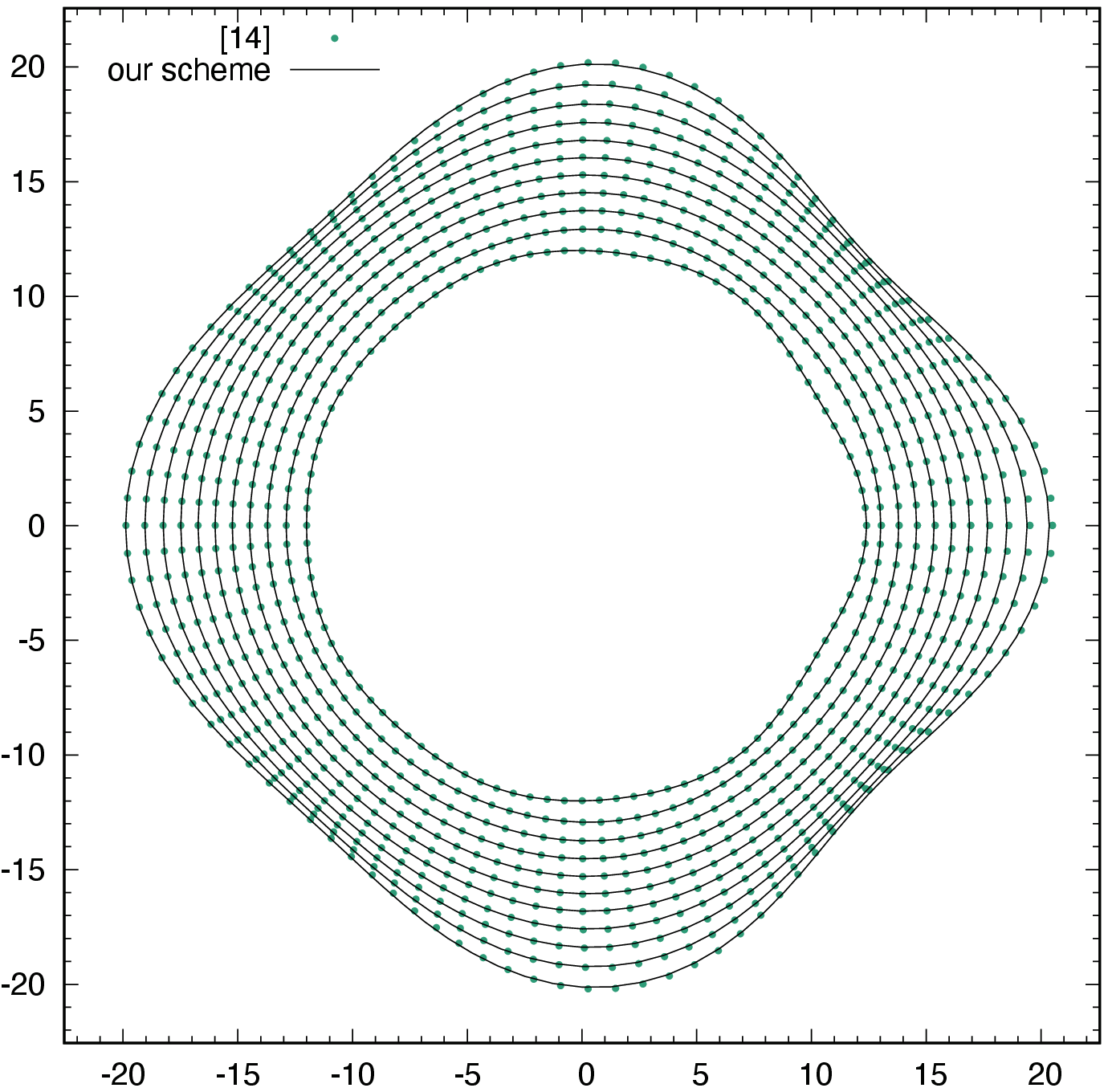}} \\ 
  \subfloat[][]{\includegraphics[width=0.35\columnwidth]{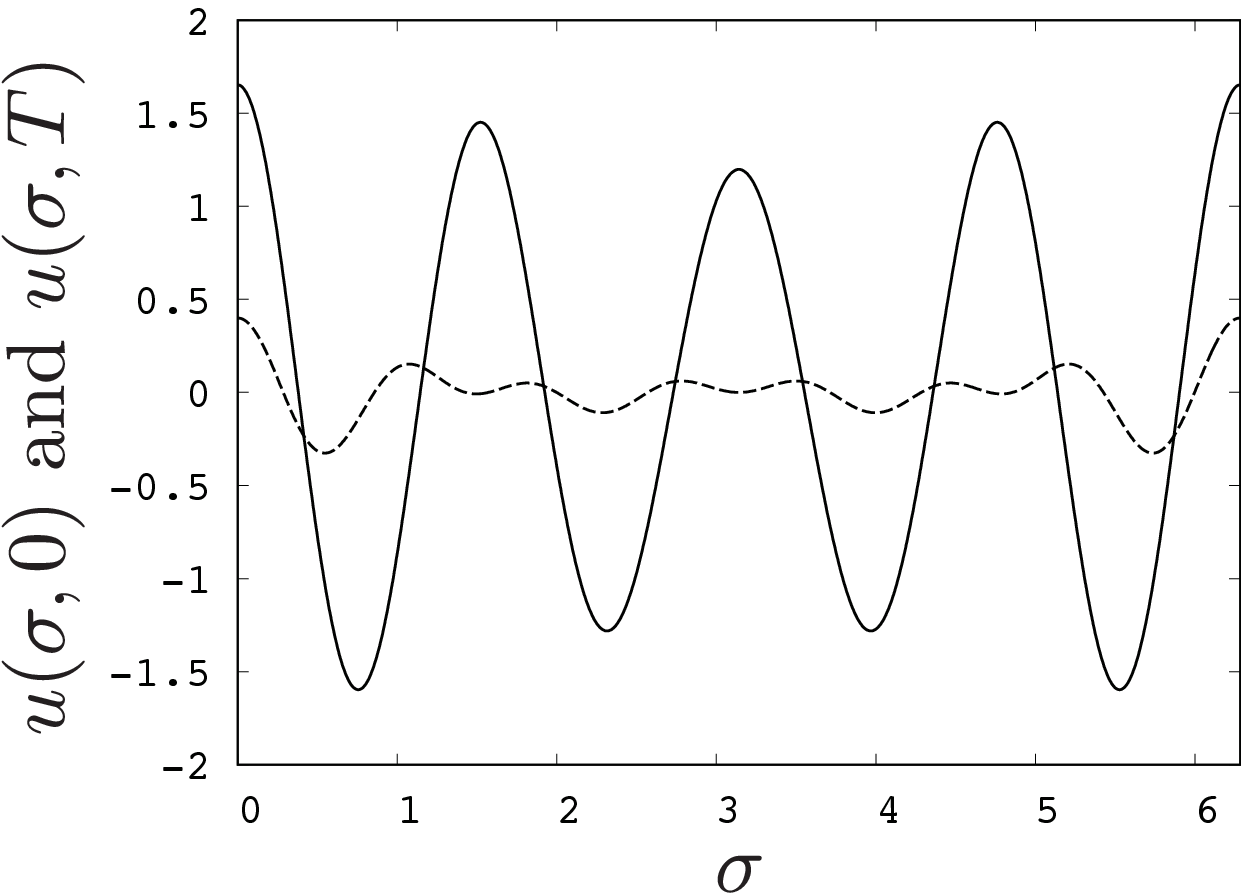}} \hspace{1cm}  
  \subfloat[][]{\includegraphics[width=0.35\columnwidth]{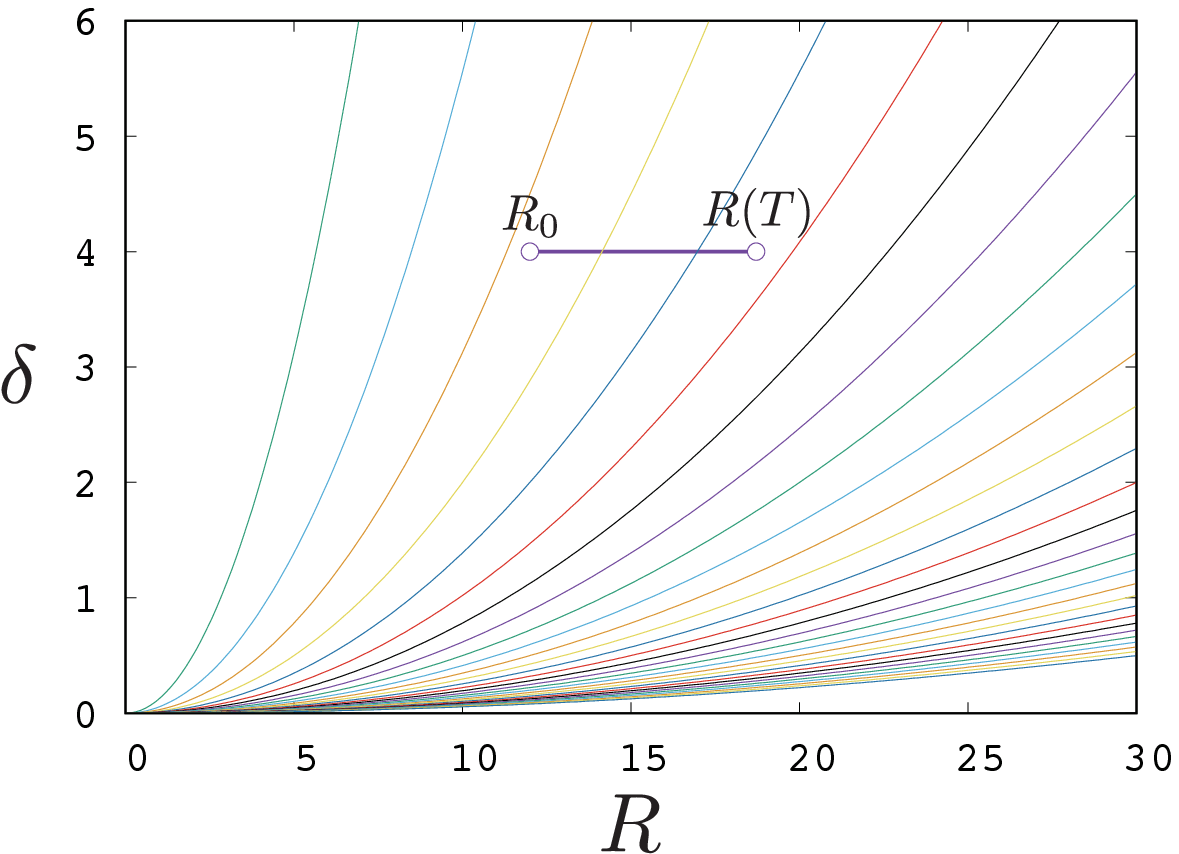}}
  \caption[]{Numerical results for the case in which $R_0 = 12.0$ and $(m_1, m_2, m_3, m_4) = (4, 5, 6, 7)$ in \eqref{eq:initial_condition}.}\label{fig:4}
\end{figure}

\begin{figure}[h]
  \centering
  \subfloat[][]{\includegraphics[width=0.48\columnwidth]{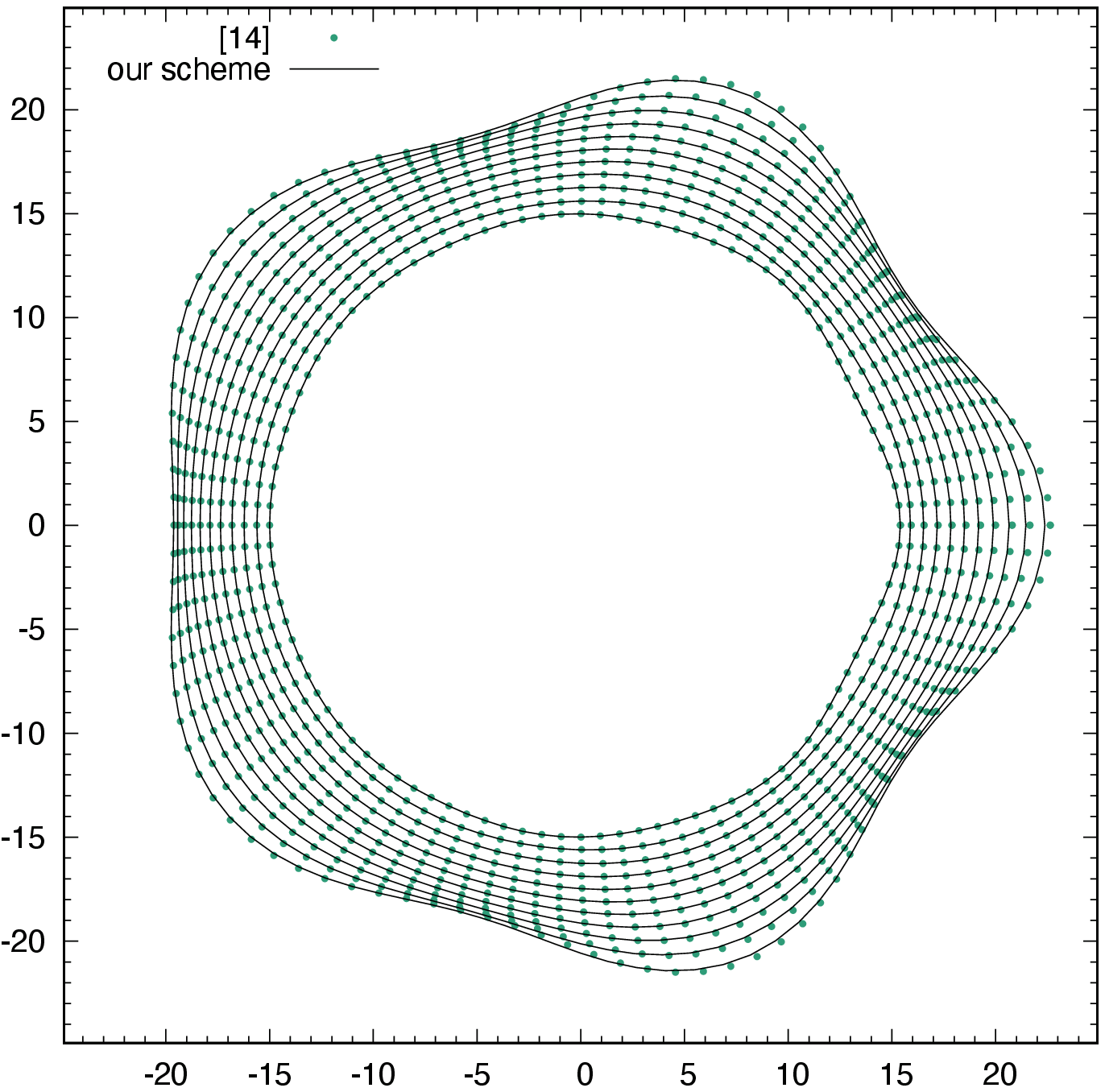}} \\ 
  \subfloat[][]{\includegraphics[width=0.35\columnwidth]{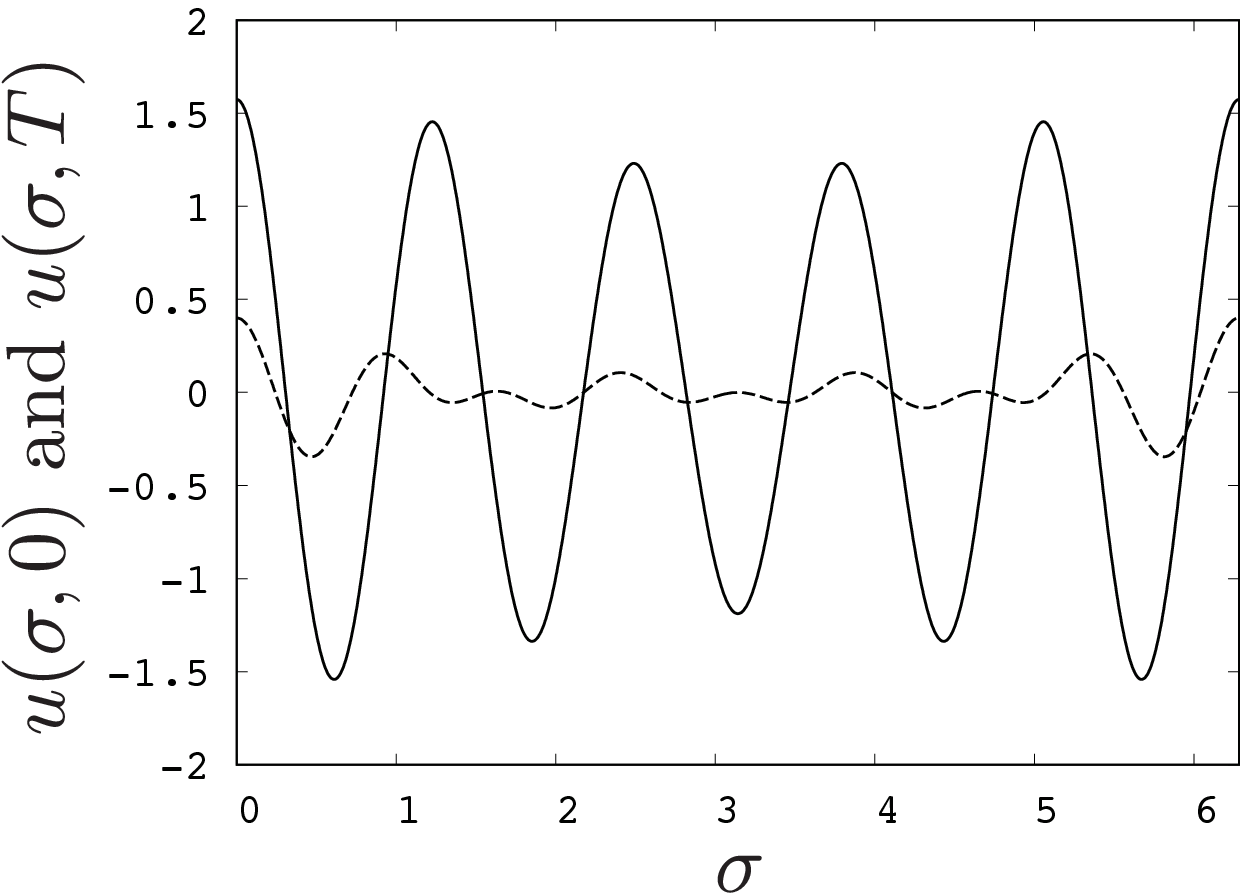}} \hspace{1cm}  
  \subfloat[][]{\includegraphics[width=0.35\columnwidth]{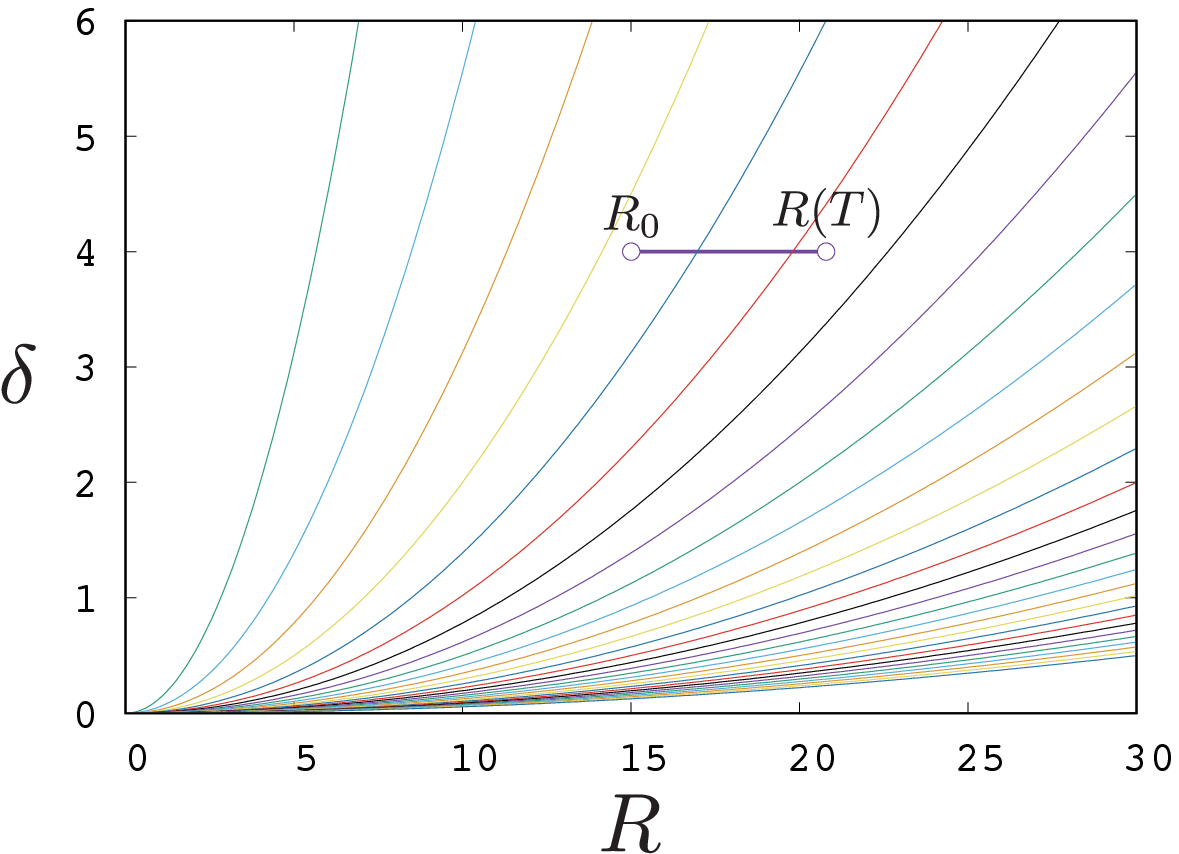}}
  \caption[]{Numerical results for the case in which $R_0 = 15.0$ and $(m_1, m_2, m_3, m_4) = (5, 6, 7, 8)$ in \eqref{eq:initial_condition}.}\label{fig:5}
\end{figure}

\begin{figure}[h]
  \centering
  \subfloat[][]{\includegraphics[width=0.48\columnwidth]{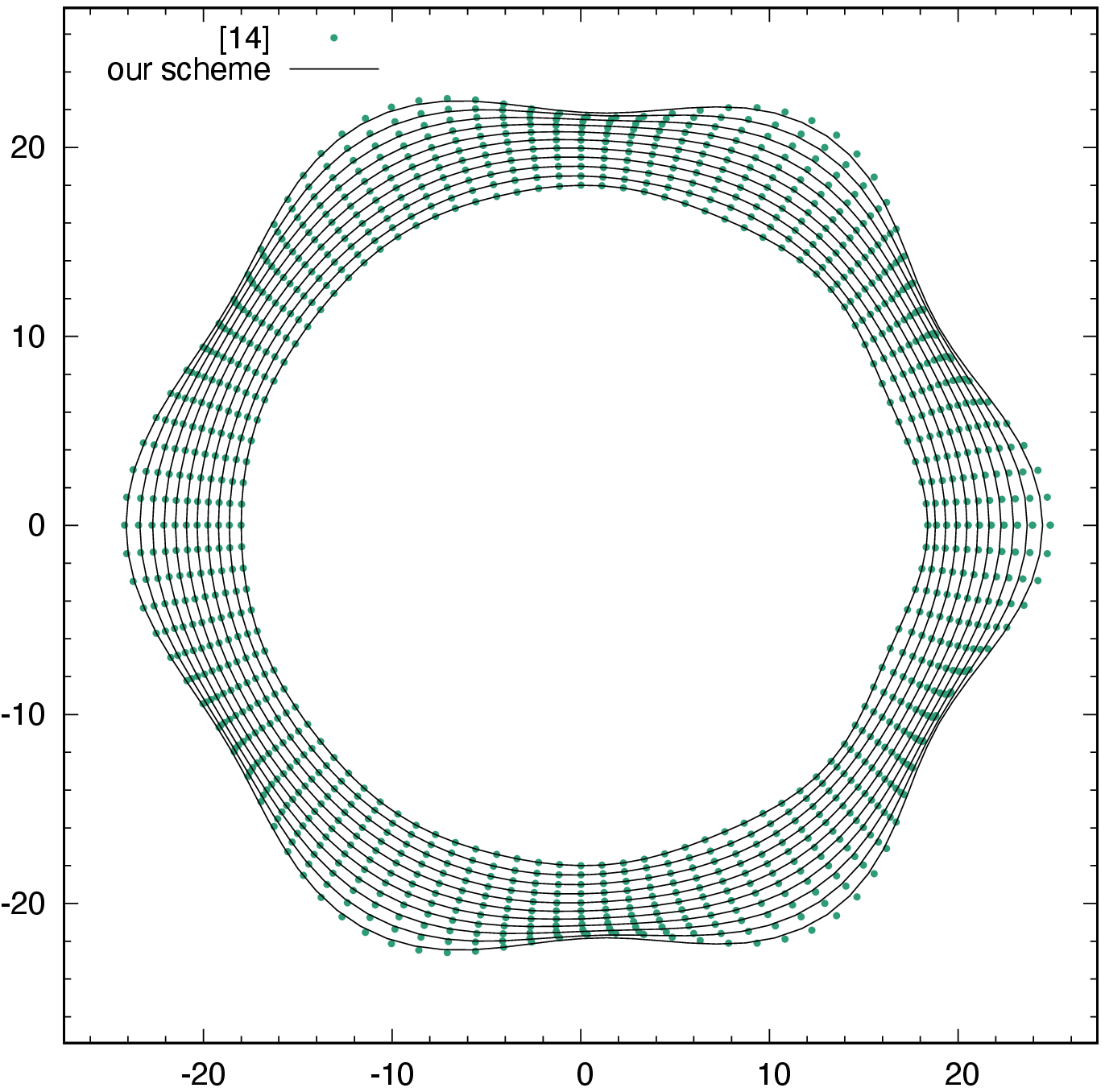}} \\
  \subfloat[][]{\includegraphics[width=0.35\columnwidth]{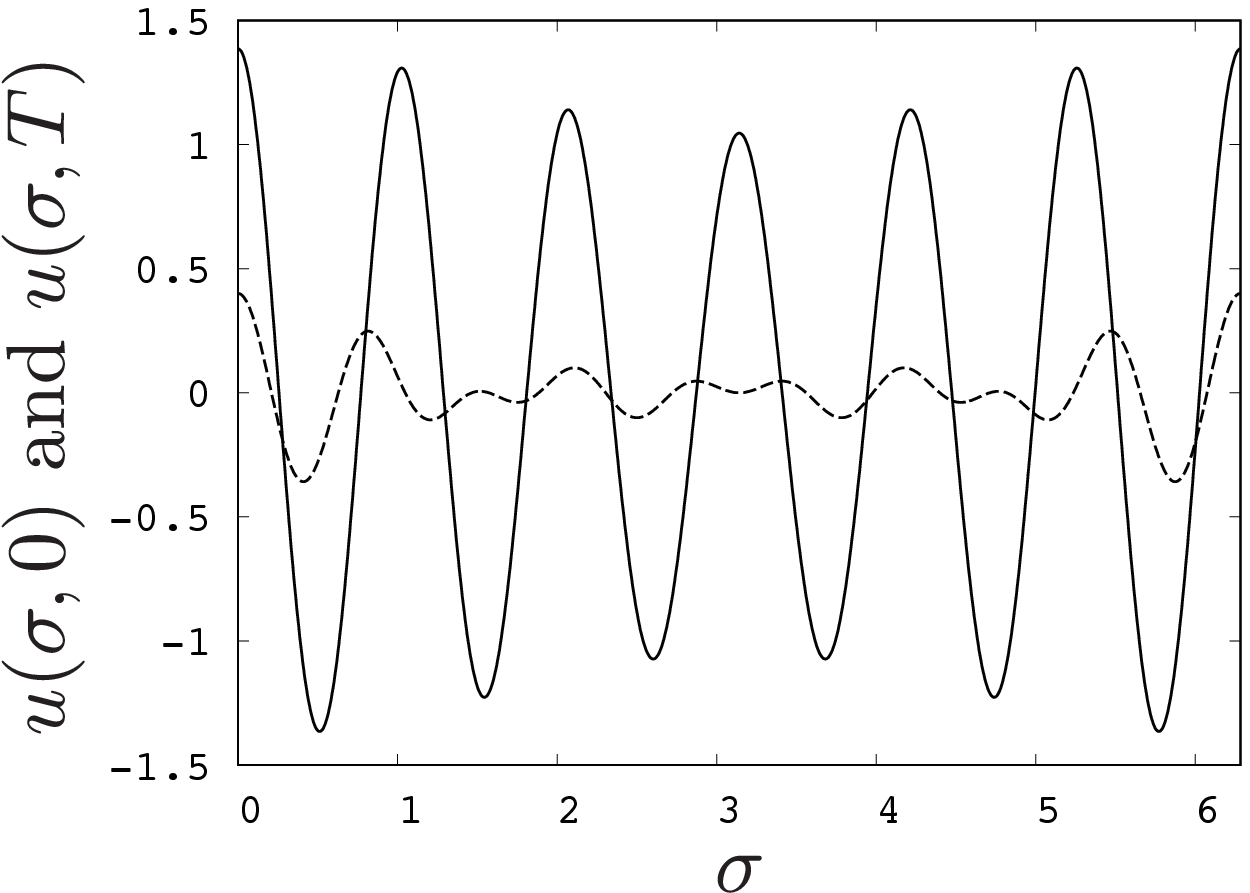}} \hspace{1cm}  
  \subfloat[][]{\includegraphics[width=0.35\columnwidth]{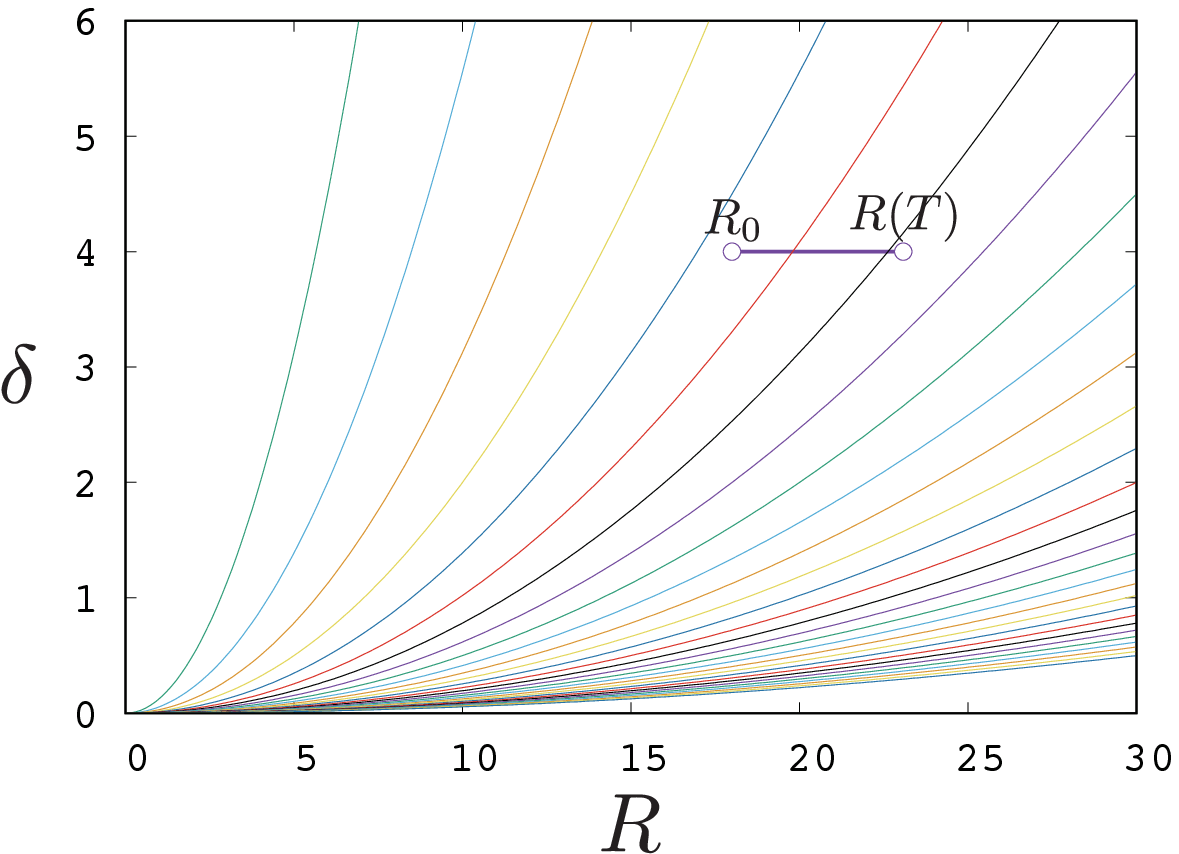}}
  \caption[]{Numerical results for the case in which $R_0 = 18.0$ and $(m_1, m_2, m_3, m_4) = (6, 7, 8, 9)$ in \eqref{eq:initial_condition}.}\label{fig:6}
\end{figure}

In this subsection, we show that the maximum wavenumber of the unstable mode can be identified a priori by appropriately choosing the parameters and initial radius.
The parameters are $dt = 0.01$, $N = 1024$, $T = 100$, $v_c = 0.001$, $\alpha = 1.5$, $\delta = 4.0$, and $p_i = 0.1$ for the initial conditions \eqref{eq:initial_condition}.
In Figs.~\ref{fig:2}--\ref{fig:6}, (a), (b), and (c) show \eqref{eq:circle} at $t = 0, 20, \dots, 100$, the numerical solution $u(\sigma, t)$, and $R(t)$ in the parameter space $(R, \delta)$ (see \cite{KUY}), respectively.
In (b) of the figures, the solid lines are the numerical solution $u(\sigma, T)$, and the dashed lines are the initial data $u(\sigma, 0)$.
The curves in (c) are called neutral stability curves $\delta = \dfrac{\alpha - 1}{m^2}R^2$ $(m \in \mathbb{N}\setminus\{1\})$, upon which the linearized operator around the trivial solution $u(\sigma, t) \equiv 0$ has at least one eigenvalue of zero.


\section{Concluding remarks}\label{sec:Conclude}

In the present paper, we proposed a simple, fast, and accurate numerical scheme for the KS equation \eqref{eq:ks_u} defined on an expanding circle. 
Our scheme is the Crank--Nicolson type finite difference scheme of \eqref{eq:ks_v} (the differential form of \eqref{eq:ks_u}), and we demonstrated the existence, uniqueness, and second-order convergence.
To the best of our knowledge, with the exception of our graph approach, there are no convergence results of the numerical scheme for the interfacial equation \eqref{eq:interfacial} in a parametric approach or a level set approach.

As mentioned above, our scheme is fast because the time increment $k$ can be taken, such as $k = o(h^{\frac{1}{4}})$, in comparison with the standard increment $k=O(h^4)$ for an explicit scheme (see Theorem~\ref{thm:2}).
In addition, our scheme is more accurate than the numerical method of \eqref{eq:interfacial} as long as the amplitude of the solution is sufficiently small because the experimental order of convergence (the so-called EOC) of the numerical method is $1$ (see \cite{KKUYB}).

Owing to the derivation of the equations \eqref{eq:ks_u}, if $\alpha$ is sufficiently small, it is expected to be a good approximation of the solution to the original interfacial equation \eqref{eq:interfacial}.
Indeed, in subsection~\ref{sec:Numerical_wavenumber}, we can see that the wavenumbers of the numerical solutions of \eqref{eq:interfacial} and \eqref{eq:ks_u} are consistent.
Therefore, we insist that our model \eqref{eq:ks_u} captures well the instability of the solution to the interfacial model.
In other words, we can determine a priori the unstable mode of the maximum wavenumber by choosing the parameters and initial radius regarding the parameter space of $(R, \delta)$ (see Fig.~\ref{fig:neutral}).

By contrast, as shown in Fig.~\ref{fig:1}, the green points (numerical solutions to \eqref{eq:interfacial}) are eventually separated from the solid lines (numerical solutions to \eqref{eq:ks_u}). 
As the reason for such separation, the self-intersection of the moving curves is allowed in the interfacial model, whereas it cannot occur in our graph model.
This observation suggests that a finite-time blow-up occurs for a solution $v(\sigma, t) = u_\sigma(\sigma, t)$.
An analysis of the blow-up phenomenon remains as one of our future studies. 


\appendix
\section{Neutral stability curves}\label{sec:apA}

In this section, we describe a linearized stability analysis of \eqref{eq:ks_u} around the trivial solution $u(\sigma, t) \equiv 0$.
Substituting the Fourier expansion $u(\sigma, t) = \sum_{m \in \mathbb{Z}} u_{m}(t) {\mathrm e}^{\sqrt{-1} m \sigma}$, $u_m(\cdot) \in \mathbb{C}$ into \eqref{eq:ks_u}, we obtain an infinite-dimensional dynamical system:
  \begin{align}\label{eq:dynamical}
    \dot u_{m}(t) = \lambda_m u_m(t) - \dfrac{v_{c}}{2R^{2}}\sum_{\substack{m_{1} + m_{2} = m \\ m_{1} m_{2} \neq 0}}m_{1}m_{2} u_{m_{1}}(t)u_{m_{2}}(t),
  \end{align}
  where $\lambda_{\pm1} \equiv 0$ and
\begin{equation*}
  \lambda_{|m| \ge 2} = -\dfrac{\delta m^{4}}{R^{4}} + \dfrac{m^{2}}{R^{2}}\left( \alpha - 1 + \dfrac{\delta}{R^{2}} \right) - \dfrac{\alpha - 1}{R^{2}}.
  \end{equation*}
Note that $u_{-m}(t) = \bar{u}_{m}(t)$ follows from $u(\cdot,\cdot) \in \mathbb{R}$.
By solving $\lambda_m = 0$ on $R$, we obtain the neutral stability curves upon which the linearized operator of \eqref{eq:ks_u} has eigenvalues of zero.
Note that $\lambda_{\pm1}\equiv0$ holds for any $R > 0$.
However, $\lambda_m < 0$ holds for each $|m| = 2,3,\cdots$ when $R < R_* = 2\sqrt{\delta/(\alpha - 1)}$.
Therefore, the circle solution is neutrally stable as indicated in the gray region in Fig.~\ref{fig:neutral}.
In the white region where $R > R_*$, the circle is unstable except at the $2$-mode neutral-stability curve because $\lambda_m > 0$ holds for at least one integer $m \in \{\pm 2, \pm 3,\dots\}$.
In particular, when $|m| = 2$, for any fixed $\delta > 0$, the value of $R_*$ is the minimum value at which the stability of the circle solution changes from neutrally stable to unstable.


\end{document}